\documentclass[11pt]{article}

\usepackage[latin1]{inputenc}
\usepackage{amsmath, amssymb, amstext, amsfonts, amsbsy, amsthm, fancyhdr, multicol,amscd}
\usepackage{graphicx}
\usepackage{amsfonts, amscd, amssymb}
\usepackage[all]{xy}

\newtheorem{Theorem}{Theorem}[section]
\newtheorem{Lemma}[Theorem]{Lemma}
\newtheorem{Corollary}[Theorem]{Corollary}
\newtheorem{Proposition}[Theorem]{Proposition}
\newtheorem{Remark}[Theorem]{Remark}
\newtheorem{Definition}[Theorem]{Definition}

\title{M\"obius solutions of the curved $n$--body problem for positive curvature}
\author{\bf Ernesto P\'erez Chavela \\
         Departamento de Matem\'aticas \\
         UAM-Iztapalapa \\
         M\'exico, D.F. MEXICO \\
         {\tt epc@xanum.uam.mx}
        \and
        \bf  J. Guadalupe Reyes-Victoria\\
         Departamento de Matem\'aticas \\
         UAM-Iztapalapa \\
         M\'exico, D.F. MEXICO \\
         {\tt revg@xanum.uam.mx}}

\begin{document}

\maketitle

\begin{abstract}
We denote by $\mathbb{M}^2_R$ a two dimensional space of constant positive Gaussian curvature.
With methods of M\"obius geometry and using the classification
 of the M\"obius group of automorphisms  ${\rm \bf  Mob}_2 \, (\widehat{\mathbb{C}})$
of the Riemman sphere $\widehat{\mathbb{C}}=\mathbb{M}_R^2 \cup \{\infty\}$,
we give algebraic conditions for the existence of M\"obius
solutions on  $\widehat{\mathbb{C}}$, getting a complete classification
of them.  We show several families of this kind of solutions.
\end{abstract}

\footnote*{MSC: Primary 70F15, Secondary 53Z05}

\footnote*{Keywords: The $n$--body problem, spaces of constant curvature, M\"obius geometry.}

\section{Introduction}
\label{sec:intro}

The formalization of Euclidian Geometry from {\it the Elements}
raised many important questions some of them about the definition of basic geometric objects.
In this way one can analyze the concept of a {\it  point} whose primary definition  was ``an object such that did not
have any parts", or as {\it the line}, which was defined as ``the object with length but without wide".
In the Nineteenth century several mathematicians found dark
these kind of definitions and show objects in suitable spaces
where they are meaningless, as the Peano spaces-filling curves.

After that, the definitions in modern
geometry are written as  mathematical objects which do not need
a primary key, because such primary definition
require of another  definition of an even more basic object
and so on.

The axiomatization of modern geometry is then
more suitable than having to define fundamental objects. In the
nineteenth century David Hilbert proposed the axiomatization of the Geometry by 21 axioms that are
fundamentals for the creation of modern non-Euclidean geometries.

In 1871, Felix Klein argued that the
classes of Euclidean and Non-euclidian objects can be studied in a
suitable projective space. In this sense the Projective Geometry establishes  a
self-reliance among the respective theories. In his work, F. Klein
proved that in order that the Euclidean geometry be consistent (without
contradictions in its postulates), is necessary and sufficient that the
non-Euclidean geometries also be consistent.

In his famous Erlangen program, Klein offers a simple definition of what
is defining a Geometry in one space, in which are not considered important
concepts of point, line, surface, etc. In his
paper proposes the idea of giving an algebraic character of
definition by using the concept of a primary group
transformations of this space (bijective applications onto
itself). That is,  the geometry of the space is defined by
properties invariant under such transformations.
These are objects that define the geometry, and
the relationships between these objects build his theory.

As one example, in order to characterize the planar Euclidean geometry, we
must define a set of rotations,
reflections and translations in the plane (isometries). The main invariants
under these applications are the points and the lines, which is expected
by the everyday experience.

With this definition, Klein did not make a distinction between
the Euclidian methods and those of the analytical algebraic geometry, which were stated
before the appearance of his program.

Among all these geometries defined in such way, is the M\"obius
geometry, which can be understand  as the study of ``Euclidean space
with a point added at infinity" endowed with a conformal metric.
That is, the setting is the compactification of the standard plane
space, and the M\"obius geometry is concerned with the group of
transformations preserving the conformal form of the metric. It is
well known that this group of conformal transformations of such
space is infinite dimensional and it is the whole set of all its
holomorphic maps (followed by a possible conjugation operation).
However, the group of conformal fractional linear M\"obius
transformations (automorphisms)  of the extended complex plane is a
complex $3$--dimensional, and its action under one-dimensional
parametric subgroups generates the conics and loxodromic curves of
such space. For more details on the development of the modern geometries and in
particular, on the theory of M\"obius geometry see \cite{Sharpe}.

We use here such group of conformal fractional M\"obius transformations in $\mathbb{C}$ for analyzing
some distinguished motions of the $n$--body problem in the positively curved and compactified complex planar space.

\smallskip

We consider $n$ bodies of masses
$m_1,\cdots,m_n$ moving on a $2$--dimensional space of constant Gaussian curvature ${\bf K}$.
 It is well known that this surface is locally characterized by the
 sign of the curvature  \cite{DoCarmo} .
\begin{enumerate}
\item  If ${\bf K}>0$, the surface is the two dimensional sphere $\mathbb{S}^2_R$ of radius
$R=1/\sqrt{{\bf K}}$ imbedded in the euclidian space $\mathbb{R}^3$, or the curved plane $\mathbb{M}^2_R$, that is, the ordinary plane $\mathbb{C}$ with the canonical complex variables $(z, \bar{z})$
endowed with the conformal metric (see \cite{Perez} for more details)
\begin{equation}\label{met-c}
   ds^2= \frac{4R^4 \,dz d\bar{z}}{(R^2 + |z|^2)^2}.
\end{equation}

\item  If ${\bf K}=0$, we recover the Euclidean space $\mathbb{R}^2$.
\item  If ${\bf K}<0$, the surface is the upper part of the hyperboloid $x^2+y^2-z^2={{\bf K}}^{-1}$
 imbedded in the three dimensional Minkowski space $\mathbb{R}^3_1$,
or the hyperbolic Poincar\'e  disk $\mathbb{D}^2_R$ with the canonical complex variables $(z, \bar{z})$ and
endowed with the conformal metric (see \cite{Diacu2} for more details)
\begin{equation}\label{metneg}
   -ds^2= \frac{4R^4 \,dz d\bar{z}}{(R^2 - |z|^2)^2}.
\end{equation}
\end{enumerate}

\smallskip

In \cite{Diac-Perez}, Florin Diacu and Ernesto P\'erez-Chavela give an analytical  definition of
the so called Eulerian and Lagrangian homographic solutions for the 3-body problem both in the
unitary two dimensional sphere  $\mathbb{S}^2$  imbedded in the euclidian space $\mathbb{R}^3$  and in the
pseudo-sphere $\mathbb{L}^2_R$ imbedded in the Minkowski space  $\mathbb{R}^3_1$.
Eulerian and Lagrangian motions are a generalization of the well known Euler and Lagrange homographic orbits
of classical celestial mechanics. In \cite{Diac-Perez}, the authors classify
all types of homographic solutions, defined as the motions where the configuration of the particles
is always on the same plane and it is similar with itself for all time. In this paper we generalize the above definition in terms of the action of M\"obius  transformations on the configuration space
for the positive curvature case, that is on  the spherical plane $\mathbb{M}^2_R$.

\smallskip

The paper is organized as follows: in section \ref{sec:equamotion} we state   the equation of motions of the problem
by using complex intrinsic coordinates as in \cite{Perez} and \cite{Diac-Perez}.
In section \ref{sec:mobius-group} we state the algebraic and geometric classification
of all M\"obius transformations in:
{\it M\"obius-elliptic, M\"obius-hyperbolic, M\"obius-parabolic} and
{\it M\"obius-loxodromic},
whose action in $\mathbb{M}^2_R$ generate the corresponding conic curves of the M\"obius geometry.

\smallskip

Starting in section \ref{elliptic}, we use the algebraic classification of the M\"obius curves
in ${\rm \bf  Mob}_2 \, (\widehat{\mathbb{C}})$ to obtain a classification of motions in the
$n$--body problem in the positive curvature case. First
we define the M\"obius-elliptic  solutions,
thorough the action of a one-dimensional parametric
subgroup of isometries (and the corresponding M\"obius subgroup of transformations)
showing that they correspond to relative equilibria.
We state the algebraic conditions that the solution must hold
in order to be one of such orbits.

\smallskip

In section \ref{hyperbolic} we define the so called M\"obius-hyperbolic solutions,
via the action of the one-dimensional parametric subgroup of M\"obius hyperbolic transformations
 showing that they correspond to some particular type of the
homothetic orbits found in \cite{Diac-Perez}. We also give the algebraic conditions on the positions
of the particles for having one of such solutions and show examples of this kind of solutions in the
curved $2$--body problem.

\smallskip

In section \ref{parabolic}
we define the M\"obius-parabolic solutions. As in the
previous cases we give the necessary and sufficient conditions which such orbits must hold. We show examples
of this kind of orbits for the curved $2$ and $3$--body problem. Until we know, this is the first time that this
kind of orbits are described explicitly, remarking the advantages of using the tools of the M\"obius geometry in the
analysis of the solutions of the curved $n$--body problem.

Finally in section \ref{loxodromic}, by using a set of one-parametric conformal transformations, we
give a geometric definition of loxodromic or  homographic solutions, which generalize those
given in \cite{Diac-Perez}. This is the largest family of motions, among all them,  we define and analyze
three distinguished classes. The first class is defined by combining the hyperbolic types obtained in section
\ref{hyperbolic} with the M\"obius-elliptic ones,
arising the class of asymptotic M\"obius loxodromic motions.
The second class is a particular type of homotethetic solutions, lying on the meridians and parametrized as a
geodesic curve;  the third class is obtained by the combined type of totally geodesic solutions with the M\"obius-elliptic
solutions, arising the class of homographic M\"obius-loxodromic motions.

\section{Equations of motion}\label{sec:equamotion}

We state in this section the equations of motion for the two-dimensional
positively curved $n$--body problem.

\smallskip

For the pair of points $z_k$ and $z_j$ in $\mathbb{M}^2_R$ we denote the geodesic distance between them
by $d(z_k,z_j)=d_{kj}$ and define (see \cite{Perez}) the cotangent relation
\begin{equation}\label{eq:cot}
\cot_R \left(\frac{d_{kj}}{R}\right)=\frac{2(z_k \bar{z}_j+z_j
\bar{z}_k)R^2+
(|z_k|^2-R^2)(|z_j|^2-R^2)}{[\Theta_{1,(k,j)} (z, \bar{z})]^{1/2}}, \\
\end{equation}
where
\begin{eqnarray} \label{eq:discot}
\Theta_{1,(k,j)} (z, \bar{z}) &= & 4 R^2 \, (z_j-z_k) \,(\bar{z}_j-\bar{z}_k) \, |R^2+ \bar{z}_j z_k|
\, |R^2+ \bar{z}_k z_j|\nonumber \\
&=& 4 R^2 \, |z_j-z_k|^2 \, |R^2+ \bar{z}_j z_k|^2. \nonumber \\
\end{eqnarray}

The singular set in  $\mathbb{M}^2_{R}$ for the  $n$--body problem
is the zero set of equation (\ref{eq:discot}):
 \[ \Theta_{1,(k,j)} (z, \bar{z})=0, \]
from here we obtain the following sets:
\begin{enumerate}
\item   The {\it singular collision set}  given by  $\Delta^{+} = \cup_{kj} \,\,  \Delta^{+}_{kj}$, where,
\begin{equation}
\label{eq:collisionset} \Delta^{+}_{kj}=\{\mathbf{z}= (z_1, z_2,\cdots, z_n) \in \mathbb{C}^n \, | \, z_k = z_j \}.
\end{equation}
\item   The {\it singular antipodal set}    given by $\Delta^{-}= \cup_{kj} \, \, \Delta^{-}_{kj}$, where,
\begin{equation}
\label{eq:antipodalset} \Delta^{-}_{kj}=\left\{\mathbf{z}= (z_1, z_2,\cdots, z_n) \in \mathbb{C}^n \, | \, z_k=
\frac{-R^2}{|z_j|^2} \, z_j \right\}.
\end{equation}
\end{enumerate}

We define the {\it total singular set} of  the problem as
\begin{equation}
\label{eq:totalsingularset} \Delta=\Delta^{+} \cup  \,  \Delta^{-}.
\end{equation}

Let  $\mathbf{z}=(z_1,z_2,\cdots,z_n) \in \mathbb{C}^n$ be the (complex) position of
$n$ point particles with masses $m_1,m_2,\cdots,m_n>0$ in the
space $\mathbb{M}^2_{R}$.  We assume that the particles are moving
under the action of the Lagrangian
\begin{equation}\label{eq:lagrangiangeral}
 L_R (\mathbf{z}, \bar{\mathbf{z}},\dot{\mathbf{z}},\dot{\bar{\mathbf{z}}})=
 K_R(\mathbf{z},\bar{\mathbf{z}},\dot{\mathbf{z}},\dot{\bar{\mathbf{z}}}) + U_R (\mathbf{z}, \bar{\mathbf{z}}),
\end{equation}
 where
 \begin{equation}\label{eq:kinetic-energy}
K_R= K_R(\mathbf{z},\bar{\mathbf{z}},\dot{\mathbf{z}},\dot{\bar{\mathbf{z}}}) = \frac{1}{2}\sum_{k=1}^n m_k
\lambda (z_k, \bar{z}_k) \, |\dot{z}_k|^2
\end{equation}
is the kinetic energy,
\begin{eqnarray}\label{eq:potesf}
U_R &=& U_R (\mathbf{z}, \bar{\mathbf{z}}) = \frac{1}{R}\sum_{1\leq k < j \leq n}^n m_k m_j
\cot_R
\left(\frac{d_{kj}}{R}\right)\\
&=& \frac{1}{R} \sum_{1\leq k < j \leq n}^n m_k m_j \frac{2(z_k
\bar{z}_j+z_j \bar{z}_k)R^2+
(|z_k|^2-R^2)(|z_j|^2-R^2)}{2 R \, |z_j-z_k| \, |R^2+ \bar{z}_j z_k|}, \nonumber
\end{eqnarray}
is the cotangent force function (i.e. the negative of the potential)
defined in the set $(\mathbb{M}^{2}_{R})^n \setminus \Delta$,  and
\begin{equation}\label{eq:conforesf}
\lambda (z_k, \bar{z}_k)= \frac{4R^4}{(R^2+|z_k|^2)^2}
\end{equation}
is the conformal function for the Riemannian metric.

\smallskip

The  solution of the corresponding Euler-Lagrange equations associated to the Lagrangian
(\ref{eq:lagrangiangeral}) satisfies the following
system of second order ordinary differential equations
\begin{equation}\label{eq:motiongral}
 m_k \ddot{z}_k -\frac{2 m_k \bar{z}_k\dot{z}_k^2}{R^2+ |z_k|^2} = \frac{2}{\lambda (z_k, \bar{z}_k)} \,
 \frac{\partial U_R }{\partial \bar{z}_k}(z, \bar{z}),
\end{equation}
where
\begin{eqnarray}\label{eq:gradmet}
 \frac{\partial U_R}{\partial \bar{z}_k} &=&
 \sum_{j=1, j \neq k}^n \frac{m_k m_j \, ( R^2+
 |z_k|^2)(|z_j|^2+R^2)^2(R^2+\bar{z}_jz_k)(z_j-z_k)}{4 R^2 \, |z_j-z_k|^3 \, |R^2+ \bar{z}_j z_k|^3} \nonumber \\
 &=&  \sum_{j=1, j \neq k}^n \frac{m_k m_j \, ( R^2+
 |z_k|^2)(|z_j|^2+R^2)^2}{4 R^2 \, |z_j-z_k| \, |R^2+ \bar{z}_j z_k| \, (\bar{z}_j-\bar{z}_k) \, (R^2+\bar{z}_k z_j)},
\end{eqnarray}
for $k=1,2, \cdots, n$.

\begin{Remark}
We observe that in equation (\ref{eq:motiongral}), the left hand side correspond to the equations of the
geodesics. This means that if the potential is constant, then the particles move along geodesics.
\end{Remark}

\section{The M\"obius group ${\rm \bf  Mob}_2 \, (\widehat{\mathbb{C}})$}\label{sec:mobius-group}
In this section we give a complete algebraic classification of all  M\"obius transformations defined on the
extended complex plane denoted by  $\widehat{\mathbb{C}}= \mathbb{M}_R^2 \cup \{\infty\}$, wich correspond to
 the rounded Riemann sphere of radius $R$ endowed with the metric (\ref{met-c}).  We start remembering the well known definition.

\smallskip

\begin{Definition}
A M\"obius transformation is a fractional linear transformation $
f_A : \widehat{\mathbb{C}}  \to \widehat{\mathbb{C}}$, defined by the rule
\[f_A (z) = \frac{a z+b}{c z + d}, \]
where $a,b,c,d \in \mathbb{C}$ and $ad-bc =1$.
\end{Definition}

It is easy to verify that the set of these automorphisms with the composition form a group denoted by
 ${\rm \bf  Mob}_2 \, (\widehat{\mathbb{C}})$.

\smallskip

Recalling that the complex special linear group of $2 \times 2$-matrix is defined by
\[ {\rm SL}(2, \mathbb{C})= \{A \in {\rm GL}(2, \mathbb{C}) \, | \, {\rm det} \, A =1 \}, \]
where ${\rm GL}(2, \mathbb{C})$ is the set of all non-singular matrices with complex entries, we
obtain that  any M\"obius transformation $f_A$ is associated to some nonsingular matrix
$A \in {\rm SL}(2, \mathbb{C})$,
\[
   A= \left(\begin{array}{cc}
    a   &  b     \\
   c &  d   \\
    \end{array}\right).
\]

Reciprocally, for any given matrix $A \in {\rm SL} (2,\mathbb{C})$,
we have a M\"obius transformation $\displaystyle f_A (z) = \frac{a z+b}{c z + d}$, such that
$\displaystyle f_A (z)= f_{-A} (z)$, which shows the isomorphism between the groups
\[ {\rm \bf  Mob}_2 \, (\widehat{\mathbb{C}}) \cong {\rm SL} (2,\mathbb{C}) \, / \{ \pm I\}.
\]

\begin{Definition}\label{def:fix-point} We say that a point $p \in \widehat{\mathbb{C}}$  is a fixed point
of the M\"obius transformation  $f_A: \widehat{\mathbb{C}} \to \widehat{\mathbb{C}}$
 if $f_A \, (p)=p$.
\end{Definition}

It is clear that the equation for a fix point is
\[ c z^2+ (d-a)z -b=0\]
and together with the relation $ad-bc=1$ we obtain the discriminant condition
\[ (a+d)^2-4= tr^2(A)  -4, \]
where $tr(A)=a+d$ denotes the trace of the matrix $A$. this means  that any  M\"obius transformation can have one or at most two fixed points.
This property is the key point to classify all elements of ${\rm \bf  Mob}_2 \, (\widehat{\mathbb{C}})$.
From here we obtain.

\smallskip

\begin{Definition}\label{def:Moebius-clasification}
 The M\"obius transformation $f_A: \widehat{\mathbb{C}} \to \widehat{\mathbb{C}}$
is called:
\begin{enumerate}

\item {\it Elliptic}, if  $0 \leq tr^2(A)<4$.

\item {\it Hyperbolic}, if  $ tr^2(A) >4$.

\item {\it Parabolic}, if $tr^2(A)= 4$.

\item {\it Loxodromic}, if $ tr^2(A) <0$ or $ tr^2(A)$ is not a real number.
\end{enumerate}

\end{Definition}

Since the determinant and the trace of one matrix are invariant under the operation of conjugation
of matrices ($D^{-1} A D$),
and since all the related matrices belong to ${\rm SL} \, (2, \mathbb{C})$, then, in order to get a classification
in ${\rm \bf  Mob}_2 \, (\widehat{\mathbb{C}})$ it is necessary and sufficient to analyze the
normal matrices of type
\[
   A= \left(\begin{array}{cc}
    \lambda   &  0   \\
   0 &  \frac{1}{\lambda}   \\
    \end{array}\right),
\]
where $\lambda$ is a simple eigenvalue of $A$, or the normal matrices of type
\[
   A= \left(\begin{array}{cc}
    1  &  b   \\
   0 &  1  \\
    \end{array}\right),
\]
for one repeated eigenvalue.

\smallskip

Now we are in conditions to state the following important result from the M\"obius geometry, which
will play a main role along this paper (you can find its proof in  \cite{Hans, Kob}).

\begin{Theorem}\label{Teo:Moebius-clasification}
 Let  $f_A: \widehat{\mathbb{C}} \to \widehat{\mathbb{C}}$ be a M\"obius transformation,
 then,
\begin{enumerate}
\item It is {\it elliptic}, if the eigenvalues satisfy $|\lambda|=1, \lambda\neq \pm 1$. In this case the
corresponding normal matrix is
\[
   \left(\begin{array}{cc}
   e^{\frac{i \theta}{2}}   &  0   \\
   0   &    e^{-\frac{i \theta}{2}}   \\
    \end{array}\right),
\]
and the associated M\"obius transformation is the rotation
\[f(z)= e^{i \theta}z, \]
around a suitable angle $\theta$.

\item It is {\it hyperbolic}, if the eigenvalues satisfy $\lambda= e^{\pm \theta}, \lambda \neq \pm 1$.
In this case the corresponding normal matrix is
\[
   \left(\begin{array}{cc}
   e^{\frac{\theta}{2}}   &  0   \\
   0   &    e^{-\frac{\theta}{2}}   \\
    \end{array}\right),
\]
and the associated M\"obius transformation is the homothetic map
\[f(z)= e^{\theta}z. \]

\item It is {\it  parabolic}, if the eigenvalues are repeated: $\lambda=1$ or $\lambda=-1$. In this case the
corresponding normal matrix is
\[
   \left(\begin{array}{cc}
    1   &  b   \\
   0   &    1   \\
    \end{array}\right),
\]
and the associated M\"obius transformation is the translation
\[f(z)= z+b. \]

\item It is {\it loxodromic}, if the eigenvalues satisfy $|\lambda^2| \neq 1$. In this case the
corresponding normal matrix is
\[
   A= \left(\begin{array}{cc}
    \lambda   &  0   \\
   0 &  \frac{1}{\lambda}   \\
    \end{array}\right).
\]
and the associated M\"obius transformation is the ``helicoidal map"
\[f(z)= \lambda^2 \, z. \]
\end{enumerate}

Moreover, the parabolic transformation leaves fix the infinity point whereas the elliptic, hyperbolic
and the loxodromic transformations leave fix both the origin of coordinates and the infinity point (The south pole
and the north pole respectively in our case).
\end{Theorem}

\begin{Remark}
We observe that any hyperbolic transformation is also a loxodromic one, actually the form of the matrices in 2
and 4 of Theorem \ref{Teo:Moebius-clasification} can be written in the same way, however we choose
the above notation to emphasize that in general in 4 the coefficient $\lambda$ could be a complex number.
\end{Remark}

As we shall show in the following sections, the action of each one-dimensional subgroup
 of the respective class of M\"obius transformations given in Theorem \ref{Teo:Moebius-clasification}
define the conic curves in the  M\"obius geometry.

\section{M\"obius-elliptic solutions}\label{elliptic}
The central part of this paper consists in to show how the algebraic classification of the  M\"obius curves in
 ${\rm \bf  Mob}_2 \, (\widehat{\mathbb{C}})$ allow us to obtain a nice classification of motions in the positive
 curvature case.
We start our analysis of the M\"obius solutions with the simplest ones, the so called M\"obius-elliptic solutions
which in fact have been widely
estudied in \cite{Diac, Diac-Perez, Perez}.

From differential geometry we know that the {\it group of proper isometries}  of
$\mathbb{M}^2_R$  is the quotient
$$ {\rm SU}(2) \, /  \, \{ \pm I \} $$
of the special unitary subgroup
\[{\rm SU}(2) = \{ A \in {\rm SL}(2,\mathbb{C}) \, | \, \, \bar{A}^T \,A= I \}  \,  \subset {\rm SL}(2, \mathbb{C}) \]
of the special linear group ${\rm SL}(2, \mathbb{C})$.
Each matrix $A \in {\rm SU}(2)$ has the form
\[
   A= \left(\begin{array}{cc}
    a        &  b     \\
   -\bar{b} & \bar{a}   \\
    \end{array}\right),
\]
with  $a,b \in \mathbb{C}$ satisfying  $|a|^2 +|b|^2 =1$.
See \cite{Dub} for more details.

\smallskip

We denote by $\{G_e(t) \}$ the
one-parametric subgroup of the Lie group $SU(2)/\{\pm I\}$, which acts
coordinate-wise  in $\mathbb{M}^{2n}_{R} \setminus \Delta$  and in
$\Delta$ leaving them invariant,
\[ G_e (t)=
\left(\begin{array}{cc}
    e^{it/2} & 0 \\
    0  & e^{-it/2} \\
    \end{array}\right),
\]
we associate it to {\it the Killing vector field}
\[\left(\begin{array}{cc}
    i/2 & 0 \\
    0  & -i/2 \\
    \end{array}\right)
\]
in the corresponding Lie algebra ${\it su} \, (2)$.
Such vector field defines the one-parametric family of acting elliptic M\"obius
transformations
\begin{equation}\label{eq:moebius-elliptic-initial}
f_{G_e(t)} (z) = e^{it} \, z,
\end{equation}
which are solutions of the complex differential equation
\begin{equation}\label{complex-eq}
\dot{z} = iz,
\end{equation}
which implies that all orbits are circular periodic orbits.

\begin{figure}[ht]
\centering
\includegraphics{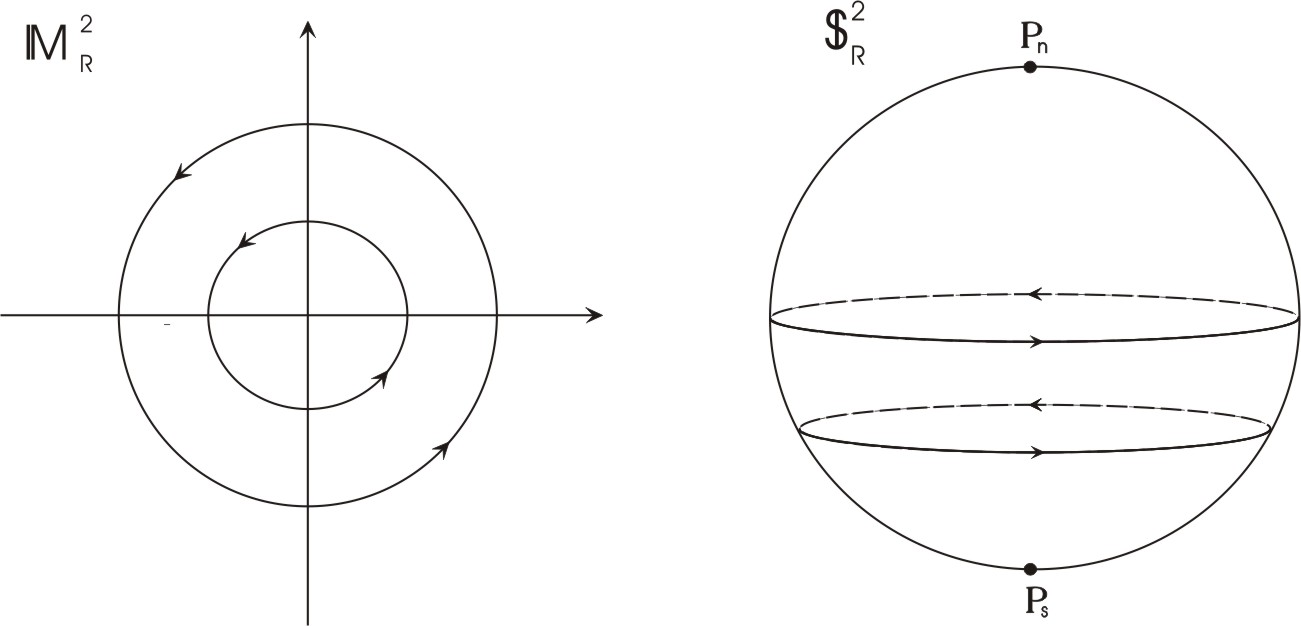}
\caption{The M\"obius-elliptic orbits on $\mathbb{M}_R^2$ and on the sphere} \label{elliptic-action}
\end{figure}

Reciprocally the above differential equation  generates the aforementioned {\it Killing vector field} associated
to the flow  $f_t(z) = e^{it}z$, and to the
one parametric subgroup of M\"obius transformations $f_{G_e(t)}$. Figure \ref{elliptic-action} shows the
M\"obius-elliptic orbits of the action of the one-parametric subgroup $\{G_e(t)\}$ in $\mathbb{M}_R^2$ and
the corresponding circular orbits in the two dimensional sphere.

\begin{Definition} A solution $\mathbf{z}(t)=(z_1(t), z_2(t), \cdots, z_n(t))$ of equation (\ref{eq:motiongral})
is called  M\"obius-elliptic if its invariant under the one parametric
subgroup of isometries (\ref{eq:moebius-elliptic-initial}).
\end{Definition}

All the above prove the following result.

\begin{Theorem} \label{match}
The M\"obius-elliptic solutions of the positive curved $n$--body problem founded with algebraic techniques
correspond with the relative equilibria, defined in \cite{Perez} as those solutions  $z(t)$ of (\ref{eq:motiongral})
which are invariant relative to the subgroup $\{G_e(t) \}$.
\end{Theorem}

Then, in order to obtain the  corresponding M\"obius-elliptic solutions in $\mathbb{M}^2_{R}$,
we must analyze just the class of M\"obius  transformations  given by
\begin{equation}\label{expo}
 w_k (t) = e^{it} z_k (t),
 \end{equation}
where $z(t) = (z_1(t), \cdots, z_n(t))$  is a solution of equation
(\ref{eq:motiongral}).

\smallskip

Introducing (\ref{expo}) in (\ref{eq:motiongral}) we can obtain the algebraic equations which characterize all
relative equilibria which, for the above Theorem they correspond to the M\"obius-elliptic solutions. This result has
been proved in \cite{Perez}. In order to have a self-contained paper, we reproduce it here without proof.

\begin{Theorem} Consider $n$ point particles
with masses $m_1,m_2, \cdots, m_n>0$ moving in $\mathbb{M}^2_{R}$. A
necessary and sufficient condition for  the solution $\mathbf{z}(t)=(z_1(t),
z_2(t), \cdots, z_n(t))$ of (\ref{eq:motiongral}) to be a M\"obius-elliptic solution (relative
equilibrium) is that the coordinates satisfy the following system
given by the rational functions:
\begin{equation} \label{eq:rationalsystem}
\frac{2 \, R^6(r_k^2-R^2) z_k }{ (R^2+ r_k^2)^4  } =
 \sum_{j=1, j \neq k}^n \frac{m_j \, (r_j^2+R^2)^2(R^2+\bar{z}_jz_k)(z_j-z_k)}{|z_j-z_k|^3 \, |R^2+ \bar{z}_j z_k|^3}
\end{equation}
where $|z_l(t)|=r_l \in [0, \pi R)$, and the velocity in each particle is given by the relation
\begin{equation}\label{eq:velocities-relative}
2i \dot{z}_k = z_k,
\end{equation}
for $k=1,2,\cdots, n$. Moreover, all obtained solutions
$z_k=z_k(t)$ are circular periodic orbits. \label{thm:existence}
\end{Theorem}

\subsection{Examples of M\"obius-elliptic solutions}\label{sec:relative-equilibria}
As we have seen in Theorem \ref{match},  the relative equilibria in the positively curved $n$--body problem
correspond to the M\"obius-elliptic solutions, so in \cite{Diac, Diac-Perez} the reader can find several families of relative equilibria defined on the sphere, in \cite{Perez} you can find must of the relative equilibria in the two and three body problems defined on $\mathbb{M}^2_{R}$. All these examples correspond in our context to
M\"obius-elliptic solutions.

\section{ M\"obius-hyperbolic solutions}\label{hyperbolic}

We start this section by defining the concept of {\it homothetic solution}.
We recall that in a Riemannian two dimensional surface, a pair of points $p_1$ and $p_2$
are conjugated if there exists a pair of different geodesics passing through them. We have that
in $\mathbb{M}^2_R$ any pair of antipodal points $\displaystyle z_k= \frac{-R^2}{|z_j|^2} \, z_j $
are conjugated and the whole space is foliated by all the geodesic curves passing through such pair of points.
The set of all such curves is called the {\it geodesic conjugated class foliation}.

\begin{Definition} A solution $\mathbf{z}(t)=(z_1(t), z_2(t), \cdots, z_n(t))$ of equation (\ref{eq:motiongral})
is called  homothetic if all the particles move on curves whose path belong to the same
geodesic conjugated class  foliation of one pair of conjugated (antipodal) points.
\end{Definition}

From the Principal Axis Theorem which states that any rotation in $\mathbb{R}^3$ is around
a fix axis,
we can assume that one point is the origin of coordinates $z=0$
with conjugated point $z=\infty$,  and the geodesic foliation of $\mathbb{M}^n_R$
in this case, is the singular set of straight lines passing
through such point, called by short {\it meridians}.

\begin{Remark}
In general  the path of one particle moving along one homothetic solution is parametrized by
$z_k(t)= \phi (t) z_{k,0}$  a suitable real function $\phi= \phi(t)$ which  holds the equations of
motion  (\ref{eq:motiongral}). We observe that with this parametrization the length of the velocity
along the curve can vary, however the  geodesics are always parameterized such that their tangent
vectors have constant speed (see \cite{DoCarmo} for more details). In other words the path of one particle 
moving along a homothetic solution not necessarily does it with constant speed, that is in a geodesic way.
\end{Remark}

Among the whole set of homothetic solutions of (\ref{eq:motiongral}) there is a subclass generated
by the action of a particular  one-parametric subgroup of hyperbolic M\"obius transformations.
Let us denote by $\{G_h(t) \}$ the
one-parametric subgroup of ${\rm \bf  Mob}_2 \, (\widehat{\mathbb{C}})$ which acts
coordinatewise  in $\mathbb{M}^{2n}_{R} \setminus \Delta$  and in
$\Delta$ leaving them invariant, and defined by
\[ G_h (t)=
\left(\begin{array}{cc}
    e^{t/2} & 0 \\
    0  & e^{-t/2} \\
    \end{array}\right),
\]
associated to the hyperbolic M\"obius group of  transformations
\begin{equation}\label{eq:Moebius-hyperbolic}
f_{G_h(t)} (z) = e^{t} \, z,
\end{equation}
also called {\it hyperbolic  group}.

\begin{figure}[ht]
\centering
\includegraphics{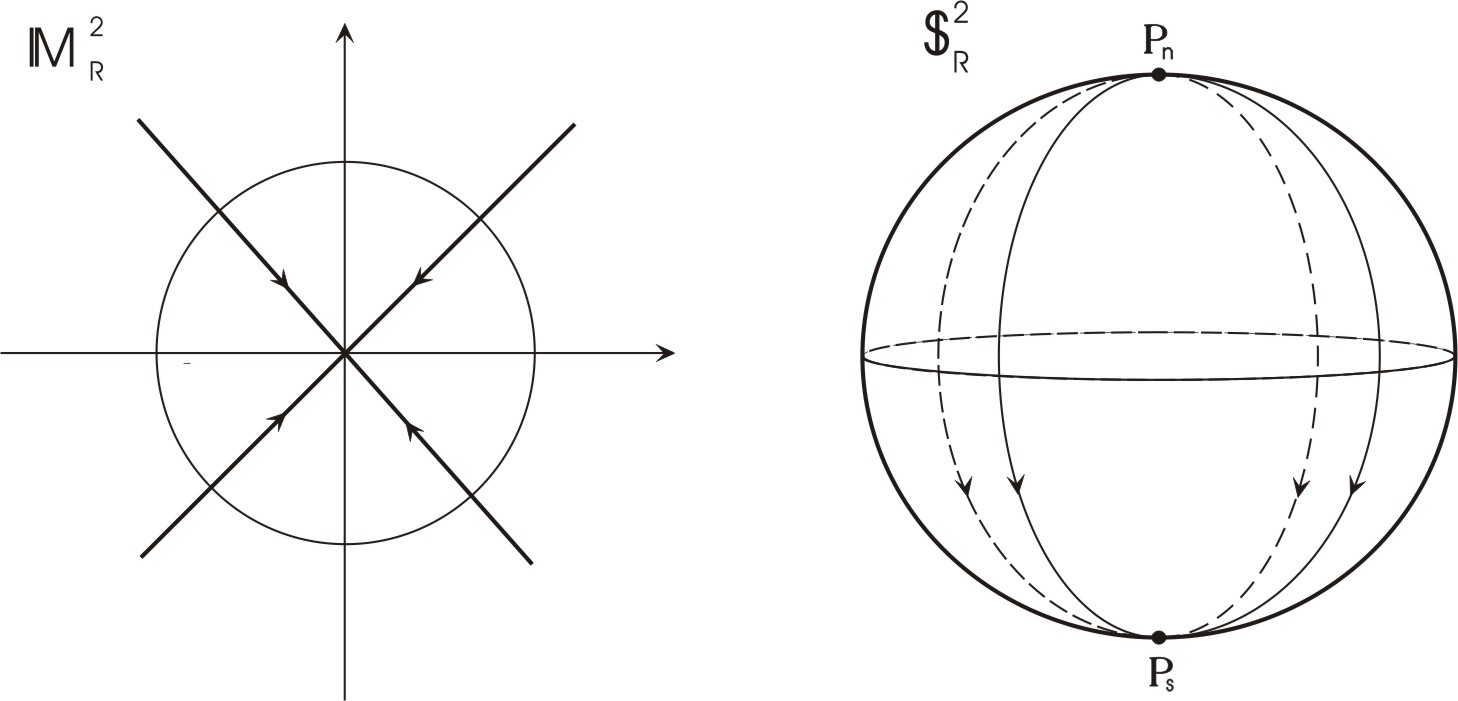}
\caption{The M\"obius-hyperbolic  orbits on $\mathbb{M}_R^2$ and on the sphere.} \label{hyperbolic-action}
\end{figure}

Figure \ref{hyperbolic-action} shows the
M\"obius-hyperbolic orbits of the action of the one parametric subgroup $\{G_h(t)\}$ in
$\mathbb{M}_R^2$ and the corresponding hyperbolic orbits in the two dimensional
sphere, which become the geodesics (great circles) of such space.

\begin{Definition} A homothetic solution $\mathbf{z}(t)=(z_1(t), z_2(t), \cdots, z_n(t))$ of equation (\ref{eq:motiongral})
is called M\"obius-hyperbolic if it is invariant under the one parametric subgroup of M\"obius
transformations (\ref{eq:Moebius-hyperbolic}).
\end{Definition}

We state the main result of this section.

\begin{Theorem} Consider $n$ point particles
with masses $m_1,m_2, \cdots, m_n>0$ moving in $\mathbb{M}^2_{R}$. A
necessary and sufficient condition for  the solution $\mathbf{z}(t)=(z_1(t),
z_2(t), \cdots, z_n(t))$ of (\ref{eq:motiongral}) to be a M\"obius-hyperbolic solution
is that the coordinates satisfy the following system
given by the rational functions:
\begin{equation} \label{eq:rationalsystem-hyp}
\frac{2 \,R^6(R^2-|z_k|^2) z_k }{(R^2+ |z_k|^2)^4} =
 \sum_{j=1, j \neq k}^n \frac{m_j \, (|z_j|^2+R^2)^2(R^2+\bar{z}_jz_k)(z_j-z_k)}{|z_j-z_k|^3 \, |R^2+ \bar{z}_j z_k|^3}
\end{equation}
and the velocity in each particle is given by the relation
\begin{equation}\label{eq:velocities-relative}
2 \dot{z}_k =- z_k,
\end{equation}
for $k=1,2,\cdots, n$. Moreover, all the solutions
$z_k=z_k(t)$ are forward asymptotic to the origin of coordinates,
which implies that there are not collisions between the particles in finite time.
\label{thm:existence-hyp}
\end{Theorem}

{\bf Proof.} If we suppose that $w_k  = e^{t} \, z_k$ is a solution of  equation (\ref{eq:motiongral}),
then by differentiating we obtain
\begin{eqnarray}
\dot{w}_k  &=& (\dot{z}_k + z_k ) \, e^{t} \nonumber \\
\ddot{w}_k  &=& (\ddot{z}_k + 2 \dot{z}_k+ z_k ) \, e^{t}, \nonumber \\
\end{eqnarray}
which when we substitute into the equation (\ref{eq:motiongral}) gives us the relation
\begin{equation}
m_k (\ddot{z}_k + 2 \dot{z}_k+ z_k ) \, e^{t} -\frac{2 m_k \bar{z}_k (\dot{z}_k + z_k )^2 \, e^{3t}}{R^2+ |z_k|^2\, e^{2t}} =
\frac{2}{\lambda (z_k \, e^{t}, \bar{z}_k\, e^{t})} \,  \frac{\partial U_R }{\partial \bar{z}_k}(z, \bar{z}) \, e^{-t},
\end{equation}
where $\lambda (z_k, \bar{z}_k)$ is the conformal function defined in (\ref{eq:conforesf}).

By evaluating at $t=0$ we obtain the condition for the infinitesimal generator of the vector field
associated to the solution $z_k=z_k(t)$ in one arbitrary point. This gives us the relation
\begin{equation}
m_k (\ddot{z}_k + 2 \dot{z}_k+ z_k ) -\frac{2 m_k \bar{z}_k (\dot{z}_k + z_k )^2 }{R^2+ |z_k|^2} =
\frac{2}{\lambda (z_k, \bar{z}_k)} \,  \frac{\partial U_R }{\partial \bar{z}_k}(z, \bar{z}),
\end{equation}
and using that $z_k=z_k(t)$ holds equation (\ref{eq:motiongral}), then this last equation becomes into
\begin{equation}\label{eq:generator:hyperbolic}
2 \dot{z}_k+ z_k =\frac{2 \bar{z}_k (2 z_k \, \dot{z}_k + z_k^2)}{R^2+ |z_k|^2}=
\frac{2  |{z}_k|^2 (2  \dot{z}_k + z_k)}{R^2+ |z_k|^2}.
\end{equation}

The equation (\ref{eq:generator:hyperbolic}) holds for the infinitesimal conditions
\begin{equation}\label{eq:generator:hyperbolic2}
|z_k|=R,
\end{equation}
which corresponds to the geodesic equator, or for
\begin{equation}\label{eq:generator:hyperbolic3}
 \dot{z}_k = -\frac{z_k}{2},
\end{equation}
which indicate that the solution through the point $z(0)$ is in the direction of the meridian passing by
the origin of coordinates and the given point.

If we derive the infinitesimal condition (\ref{eq:generator:hyperbolic3}),
we obtain $\displaystyle \ddot{z}_k = \frac{z_k}{4}$, which when is substituted in equation of motion
(\ref{eq:motiongral}) allows to the system (\ref{eq:rationalsystem-hyp}),  this ends the proof of the Theorem. \qed


\subsection{Examples of  M\"obius-hyperbolic solutions}\label{sec:hyperbolic-solutions}
For the two body problem we have the following result

\begin{Theorem} \label{theo:two-bodies-hyp}
Suppose that two particles with masses $m_1$ and $m_2$ and positions
$z_1(t)$ and $z_2(t)$  in $\mathbb{M}^2_{R}$
move as a M\"obius-hyperbolic solution. Then the masses are equal  if and only if  $z_1(t)=-z_2(t)$.
\end{Theorem}

{\bf Proof.} For two bodies  in $\mathbb{M}^2_{R}$ with masses $m_1$ and $m_2$
moving as a M\"obius-hyperbolic solution, the system
(\ref{eq:rationalsystem-hyp}) becomes into,
\begin{eqnarray} \label{eq:rationalsystem-hyp-two-bodies}
\frac{2 \,R^6(R^2-|z_1|^2) z_1 }{(R^2+ |z_1|^2)^4} &=&
 \frac{m_2 \, (|z_2|^2+R^2)^2(R^2+\bar{z}_2z_1)(z_2-z_1)}{|z_2-z_1|^3
 \, |R^2+ \bar{z}_2 z_1|^3}, \nonumber \\
 \frac{2 \,R^6(R^2-|z_2|^2) z_2}{(R^2+ |z_2|^2)^4} &=&
  \frac{m_1 \, (|z_1|^2+R^2)^2(R^2+\bar{z}_1 z_2)(z_1-z_2)}{|z_1-z_2|^3
 \, |R^2+ \bar{z}_1 z_2|^3}, \nonumber \\
\end{eqnarray}
and if we divide term to term both sides of this system avoiding singularities and conjugated points,
after a straightforward algebraic manipulation, we obtain the relation
\begin{eqnarray}\label{eq:hyp-two-bodies}
0 &=& m_1 \,z_1 (R^2-|z_1|^2) (|z_2|^2+R^2)^2 (R^2+\bar{z}_1 z_2) \nonumber \\
 &+& m_2 \,z_2 (R^2-|z_2|^2) (|z_1|^2+R^2)^2(R^2+\bar{z}_2 z_1)=0. \nonumber \\
\end{eqnarray}

If we put $z_1(t)=-z_2(t)$ into  equation (\ref{eq:hyp-two-bodies}), then necessarily the masses
must be equal.

On the other hand, if $m_1=m_2$, without lose of generality (up a rotation) we take
$z_1=r$ as one real number, then equation (\ref{eq:hyp-two-bodies}) implies that necessarily
$z_2$ is also a real number, say $z_2= \alpha$. Therefore such equation becomes into the
algebraic equation in the unknown $\alpha$,
\begin{equation}\label{eq:hyp-two-bodies-alpha}
0 =[r (R^2-r^2) (\alpha^2+R^2)^2 + \alpha (R^2-\alpha^2) (r^2+R^2)^2](R^2+r \alpha).
\end{equation}

Since we avoid conjugated points, then $R^2+r \alpha \neq 0$, and therefore we obtain
the real equation
\begin{equation}\label{eq:hyp-two-bodies-alpha2}
0 =r (R^2-r^2) (\alpha^2+R^2)^2 + \alpha (R^2-\alpha^2) (r^2+R^2)^2,
\end{equation}
which has the only real roots $\alpha = \pm r$. This ends the proof. \qed

\begin{Remark} We observe that if we put two equal masses on initial positions $z_1(t_0)=-z_2(t_0)$, then we can easily generate a M\"obius-hyperbolic solution. \end{Remark}

In \cite{Diac-Perez}, the authors show that a necessary and sufficient condition for
having a homothetic solution in the curved $3$--body problem is that the configuration be
always an equilateral triangle and that the masses be equal (they called it {\it Lagrangian} homothetic solution).
So in order to study this kind of motion  it is enough to analyze the case of equal masses. Proceeding as in the
above reference, we get the same type of
M\"obius hyperbolic solution, which also call them  {\it Lagrangian}, we omit here the details.

In \cite{Perez}, the authors prove that a necessary and sufficient condition in order to have a
collinear M\"obius-hyperbolic solution in $\mathbb{M}^2_{R}$, with one particle at the origin $z_3(t)=0$,
is that the masses of the other two particles with positions $z_1(t)$ and $z_2(t)$  be equal,
and that $z_1(t)=-z_2(t)$. This is another example of a  M\"obius-hyperbolic solution.

\section{M\"obius-parabolic  solutions}\label{parabolic}
\begin{figure}[ht]
\centering
\includegraphics{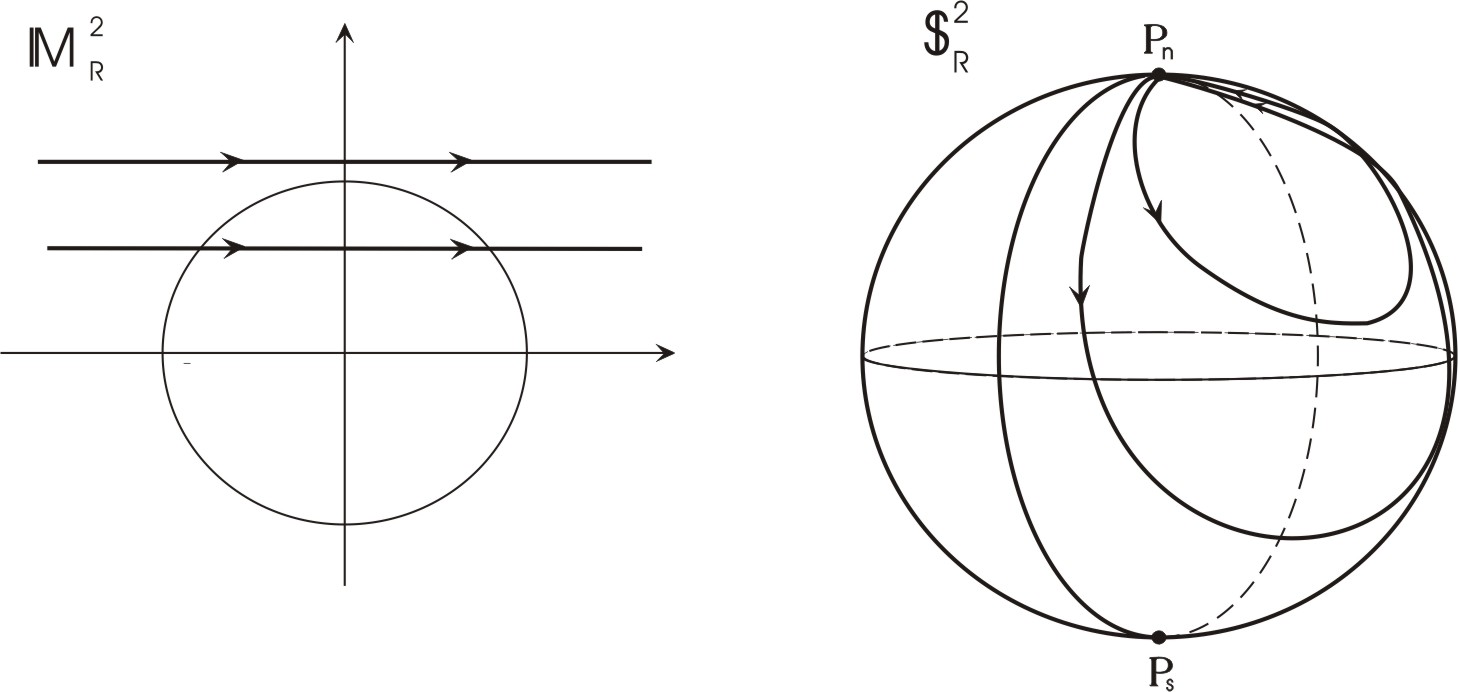}
\caption{The M\"obius-parabolic orbits on $\mathbb{M}_R^2$ and on the sphere.}
\label{parabolic-action}
\end{figure}

In this section we do the analysis  of the
 M\"obius-parabolic solutions corresponding to the solutions associated to the
one parametric subgroup
\[ G_p (t)=
\left(\begin{array}{cc}
    1 & t \\
    0  & 1  \\
    \end{array}\right),
\]
which defines the one-parametric subgroup of acting M\"obius transformations
\begin{equation} \label{eq:Moebius-parabolic}
f_{G_p (t)} (z) = z + t,
\end{equation}
in $\mathbb{M}^2_R$. Figure \ref{parabolic-action} shows the
straight lines or parabolic orbits of the action of the one parametric subgroup $\{G_p (t)\}$ in $\mathbb{M}_R^2$
and the corresponding looped M\"obius-parabolic orbits on the two dimensional sphere.

Let $\mathbf{w} =(w_1,\dots, w_n)$, with $w_k (t) = z_k (t) + t$, be the action orbit for $\mathbf{z}=(z_1,\dots, z_n)$,
which is a solution  of the equations of motion (\ref{eq:motiongral}). Then
$$
\dot{w}_k = \dot{z}_k+1 \ \ {\rm and}\ \
\ddot{w}_k = \ddot{z}_k, \ k=1,\dots, n,
$$
therefore $w$ is  a solution of system (\ref{eq:motiongral}), if and only if
$$
 m_k \ddot{z}_k  = \frac{2 m_k (\bar{z}_k +t)(\dot{z}_k +1)^2}{ R^2+ |z_k+t|^2} +
 \frac{ ( R^2+ |z_k+t|^2)^2}{2 R^4} \, \frac{\partial U_R}{\partial \bar{z}_k}, \ k=1,\dots,n,
$$
to get the above expression we have used the fact that $\displaystyle \frac{d \bar{w}_k}{d \bar{z}_k}=1$.

Once again, if we take $t=0$ for finding the condition for the infinitesimal generator
of the vector field associated to such motions we have the system of algebraic equations
$$
 m_k \ddot{z}_k  = \frac{2 m_k \, \bar{z}_k \, (\dot{z}_k +1)^2}{ R^2+ |z_k|^2} +
 \frac{ (R^2+ |z_k|^2)^2}{2 R^4} \, \frac{\partial U_R}{\partial \bar{z}_k}, \ k=1,\dots,n.
$$

Now, since $z$ is a solution of system (\ref{eq:motiongral}), we obtain the
condition for the  infinitesimal  generator of the vector field
\begin{equation} \label{eq:condiklein2}
2\dot{z}_k = -1,\ k=1,\dots,n,
\end{equation}
which holds if and only if
\begin{equation}\label{eq:principalcondklein2}
   z_k (t) = - \frac{t}{2}+ z_k(0),\ k=1,\dots,n,
\end{equation}
where $z_k(0),\ k=1,\dots, n,$ are initial conditions.

Consequently, a necessary condition for the particles $m_1,\dots, m_n$ to form a
M\"obius-parabolic solution is that they move along horizontal straight lines in
$\mathbb{M}^2_R$ with constant negative velocity, passing through the initial
conditions $z_k(0), \ k=1,\dots, n$.

\smallskip
 
We can now state the following result, whose proof follows by straightforward computations.

\begin{Theorem}\label{thm:existenceklein2} Consider $n\ge 2$ point particles
of masses $m_1,\dots, m_n>0$ moving in $\mathbb{M}^2_{R}$. Then a necessary and sufficient  condition for
the function $\mathbf{z}=(z_1,\dots, z_n)$ to be a M\"obius-parabolic solution of system
(\ref{eq:motiongral})  is that the coordinate functions satisfy the equations
\begin{equation} \label{eq:condrationalsystem-parabolic}
-\frac{4 R^6 \bar{z}_{k}}{(R^2+ |z_k|^2)^4} = \sum_{\substack{j=1\\ j\ne k}}^n
\frac{m_j \, (R^2+ |z_j|^2)^2(R^2+ z_k \bar{z}_j)(z_j-z_k)}{|z_j-z_k|^3 \, |R^2+ \bar{z}_j z_k|^3},
\end{equation}
for $k=1,\dots, n$.
\end{Theorem}

Notice that equations (\ref{eq:condiklein2}) give the conditions for the velocities of the particles
in case they form a M\"obius-parabolic solution.

In the next subsection,  using the
tools of the M\"obius geometry, we show how the M\"obius-parabolic solutions
 appear in a natural way when the one-parametric group
$\{G_p (t)\}$ acts in $\mathbb{M}_R^2$. These kind of parabolic
motions obtained with the M\"obius geometric tools could  hardly be obtained using the traditional approach.

\subsection{Examples of M\"obius-parabolic solutions}\label{subsec:parabolic-n-2}
We start this subsection by giving a concrete example of a  M\"obius-parabolic solution in the following lemma.

\begin{Lemma}\label{coro:two-bodies-parabolic}
For small $m>0$ there exist a pair of functions $z_1(t)=-z_2(t)$, such that
conform a M\"obius-parabolic solution of the problem.
\end{Lemma}

{\bf Proof.} If $\displaystyle z_1(t)=-\frac{t}{2}+\alpha R i$ and $\displaystyle z_1(t)=\frac{t}{2}-\alpha R i$
are substituted in the equations of motion (\ref{eq:motiongral}), then the condition $z_1(t)=-z_2(t)$ 
with equal masses $m_1=m_2=m$ in the particles carries such system into the single equation
\begin{equation} \label{eq:condrationalsystem-parabolic-2-b}
\frac{\bar{z}_{1}(t)}{(R^2+ r(t)^2)^6} =
- \frac{m \, z_1 (t)}{16 R^6 r(t)^3 (R^2-r(t)^2)^2}.
\end{equation}

Now taking $t=0$ in the above equation we have $\displaystyle z_1= \alpha R i$ and $r(0)=\alpha R$,from here we get the condition for $\alpha$ given by the equation
\begin{equation} \label{eq:condrational-existence-parabolic-2}
16 \alpha^3 (1-\alpha^2)^2 = \frac{ m}{R} (1+\alpha^2)^6.
\end{equation}

If we do $m=0$ (or $R=\infty$, the planar case) in equation (\ref{eq:condrational-existence-parabolic-2}),
the roots for the above equation are $\alpha=0, \pm 1$, and the particles are located at the origin of
coordinates, or on the geodesic circle.
Therefore, by using an argument of continuity  and the generic property of transversality for the smooth real 
functions on the  left and right hand sides of (\ref{eq:condrational-existence-parabolic-2}),
for fixed $R>0$ and small enough $m>0$ (or for a fixed $m>0$ and big enough $R>0$), 
it follows that there exists a solution of  such equation
depending of $m$ and $R$, say $\alpha= \alpha(R,m)$ such that $1< \alpha(R,m)<\pi$.
We define the functions $\displaystyle z_1(t)=-\frac{t}{2}+\alpha(R,m) \, R i$ and $\displaystyle z_1(t)=\frac{t}{2}-\alpha(R,m) \, R i$,
which by construction conform a M\"obius-parabolic solution. \qed

\medskip

Once we have show an example of a M\"obius-parabolic solution, we can state a general result about these kind of motions.

\begin{Theorem} \label{theo:two-bodies-parabolic}
Suppose that two particles with masses $m_1$ and $m_2$ in $\mathbb{M}^2_{R}$ located
in positions $z_1(t)$ and $z_2(t)$, with $|z_1(t)| =|z_2(t)|=r(t)$
move as  M\"obius-parabolic solutions (\ref{eq:principalcondklein2}).
Then, for $r(t)>R$, the masses are equal, if and only if,  $z_1(t)=-z_2(t)$.
\end{Theorem}

\smallskip

{\bf Proof.} For such two bodies with masses moving as a M\"obius-parabolic solution,
the system (\ref{eq:condrationalsystem-parabolic}) becomes into,
\begin{eqnarray} \label{eq:condrationalsystem-parabolic-2}
-\frac{4 R^6 \bar{z}_{1}}{(R^2+ r(t)^2)^6} &=&
\frac{m_2 \, (R^2+ z_1 \bar{z}_2)(z_2-z_1)}{|z_2-z_1|^3 \, |R^2+ \bar{z}_2 z_1|^3}, \nonumber \\
-\frac{4 R^6 \bar{z}_{2}}{(R^2+ r(t)^2)^6} &=&
\frac{m_1 \, (R^2+ z_2 \bar{z}_1)(z_1-z_2)}{|z_1-z_2|^3 \, |R^2+ \bar{z}_1 z_2|^3}, \nonumber \\
\end{eqnarray}
Since we avoid collisions,  by dividing the corresponding hand sides of first equation with those
of the second one, we obtain,
\begin{equation} \label{eq:condrationalsystem-parabolic-2-a}
\frac{\bar{z}_{1}}{\bar{z}_{2}} = - \frac{m_2 \, (R^2+ z_1 \bar{z}_2)}{m_1 \, (R^2+ z_2 \bar{z}_1)}.
\end{equation}

If we assume that $z_1(t)=-z_2(t)$ then $R^2+ z_1 \bar{z}_2=R^2+ z_2 \bar{z}_1=R^2- r(t)^2$,
and therefore equation (\ref{eq:condrationalsystem-parabolic-2-a}) implies that the masses must be equal.

Conversely, we suppose that the masses are equal and the functions $\displaystyle z_1(t)=-\frac{t}{2}+\alpha R i$
and $\displaystyle z_2(t)= \pm \frac{t}{2}+(C+\beta R i)$ conform a M\"obius parabolic solution for suitable $1<\alpha<\pi $.
By taking $t=0$ in equations (\ref{eq:principalcondklein2})
for the corresponding solutions, we consider the first particle sited on the imaginary
axis $\displaystyle z_1= \alpha R i$ and the other one, say, in the position $\displaystyle z_2(t)=C+\beta R i$ for some real numbers
$C$ and $\beta$. The condition  $\displaystyle \alpha^2 R^2= C^2+\beta^2 R^2$ on the given initial positions,
together with the condition (\ref{eq:condrationalsystem-parabolic-2-a}) implies that necessarily
$C=0$ and $\beta=\pm \alpha$. This is, both particles must initially be located on the imaginary axis.
Using again the condition (\ref{eq:condrationalsystem-parabolic-2-a}), it follows that necessarily  $z_1(t)=-z_2(t)$ for
all time $t \in (-\pi R, \pi R)$. \qed

\begin{Remark}  The above solutions in Corollary \ref{coro:two-bodies-parabolic} 
correspond on the two dimensional sphere to motions of
particles moving symmetrically on parabolic curves, leaving the north pole (as $t \to -\pi R$)
and going forward to the same point (as $t \to \pi R$) and the collision between the particles
does in finite time.
\end{Remark}

Now we shall show the existence of  M\"obius-parabolic motions
in the curved three body problem. Such solutions correspond on the two dimensional sphere
to motions of particles moving symmetrically on parabolic curves,
leaving the north pole (as $t \to -\pi R$)  and converging to the same point (as $t \to \pi R$).
The third arbitrary mass  is located always in the conjugated point, that is, at the south pole.

\begin{Lemma} \label{coro:three-bodies-parabolic}
 If $m_1=m_2=m$, then, for small $m>4 M$ there exist a pair of functions $z_1(t)=-z_2(t)$, such that
conform a M\"obius-parabolic solution of the problem, and the antipodal particles collide in the north pole in
finite time.
\end{Lemma}

{\bf Proof.} If the functions  $\displaystyle z_1(t)=-\frac{t}{2}+\alpha R i$ and 
$\displaystyle z_1(t)=\frac{t}{2}-\alpha R i$
are substituted into the equations of motion (\ref{eq:motiongral}), then the condition $z_1(t)=-z_2(t)$ with equal masses $m_1=m_2=m$
in the particles carries such system in the single equation
\begin{equation} \label{eq:condrationalsystem-parabolic-3-a}
\frac{4 R^6 \bar{z}_{1}(t)}{(R^2+ r(t)^2)^4} = - \frac{m (R^2+ r(t)^2)^2\, z_1(t)}{4 R^2 r(t)^3 (R^2-r(t)^2)^2}
+ \frac{M \, z_1(t)}{r(t)^3}.
\end{equation}

If we take $t=0$, then $\displaystyle z_1= \alpha R i$ and $r(0)=\alpha R$,  getting the condition for $\alpha$
 given by the equation
\begin{equation} \label{eq:condrational-existence-parabolic-3}
16 \alpha^3 (1-\alpha^2)^2 = \frac{(1+\alpha^2)^4 }{R}  [m \, (1+\alpha^2)^2- 4 M \,(1-\alpha^2)^2].
\end{equation}

If we do $m=0$ and $M=0$ in equation (\ref{eq:condrational-existence-parabolic-3}),
the roots for the corresponding equation are $\alpha=0, \pm 1$, and the particles are located at the origin of
coordinates, or on the geodesic circle.
By using again an argument of continuity and the generic property of transversality for 
the smooth functions on the corresponding left and right hand sides of (\ref{eq:condrational-existence-parabolic-3}),
for fixed $R>0$ and small $m>4M$, it follows that there exists a solution of the such equation
depending of $m$, $M$ and $R$, say $\alpha= \alpha(R,m,M)$ such that $1< \alpha(R,m,M)<\pi$.
Once again, we define the functions $\displaystyle z_1(t)=-\frac{t}{2}+\alpha(R,m,M) \, R i$ and
$\displaystyle z_1(t)=\frac{t}{2}-\alpha(R,m,M) \, R i$,
which by construction conform a M\"obius-parabolic solution where the antipodal particles collide in the north pole
in finite time. \qed

\begin{Theorem} \label{theo:three-bodies-parabolic}
Suppose that two particles with masses $m_1$ and $m_2$ in $\mathbb{M}^2_{R}$ located
in positions $z_1(t)$ and $z_2(t)$, with $|z_1(t)| =|z_2(t)|=r(t)$, and the third
mass $M$ located at  the origin of
coordinates, move as  M\"obius-parabolic solutions (\ref{eq:principalcondklein2}).
Then, for $r(t)>R$, the masses $m_1$ and $m_2$ are equal, if and only if,  $z_1(t)=-z_2(t)$.
\end{Theorem}

{\bf Proof.} For such two bodies with masses moving as a M\"obius-parabolic solution,
 system (\ref{eq:condrationalsystem-parabolic}) becomes into,
\begin{eqnarray} \label{eq:condrationalsystem-parabolic-3}
-\frac{4 R^6 \bar{z}_{1}}{(R^2+ r(t)^2)^4} &=&
\frac{m_2 \, (R^2+ r(t)^2)^2(R^2+ z_1 \bar{z}_2)(z_2-z_1)}{|z_2-z_1|^3 \, |R^2+ \bar{z}_2 z_1|^3}
-\frac{m_3 \, z_1}{|z_1|^3},\nonumber \\
-\frac{4 R^6 \bar{z}_{2}}{(R^2+ r(t)^2)^4} &=&
\frac{m_1 \, (R^2+ r(t)^2)^2(R^2+ z_2 \bar{z}_1)(z_1-z_2)}{|z_1-z_2|^3 \, |R^2+ \bar{z}_1 z_2|^3}
-\frac{m_3 \, z_2}{|z_2|^3}, \nonumber \\
0 &=& \frac{m_1 \, (R^2+ r(t)^2)^2  z_1}{|z_1|^3 \, R^4}
+ \frac{m_2 \, (R^2+ r(t)^2)^2 z_2}{|z_2|^3 \, R^4}, \nonumber \\
\end{eqnarray}
From third equation of (\ref{eq:condrationalsystem-parabolic-3}), we obtain the condition for the
value and position of the
first mass,
getting \begin{equation} \label{eq:condrationalsystem-parabolic-3-a}
0= m_1 z_1+ m_2 z_2,
\end{equation}
which, together with the relation $|z_1(t)| =|z_2(t)|=r(t)$ implies the assertion. \qed

\section{M\"obius-loxodromic  solutions}\label{loxodromic}

We study here two particular classes of M\"obius loxodromic motions
of the positively curved $n$--body problem.
We recall that a loxodromic transformation has associated the
matrix, which in its normal form is given by
\[
   A= \left(\begin{array}{cc}
    \lambda   &  0   \\
   0 &  \frac{1}{\lambda}   \\
    \end{array}\right),
\]
its corresponding M\"obius transformation is the helicoidal map
\[f(z)= \lambda^2 z. \]

\begin{figure}[ht]
\centering
\includegraphics{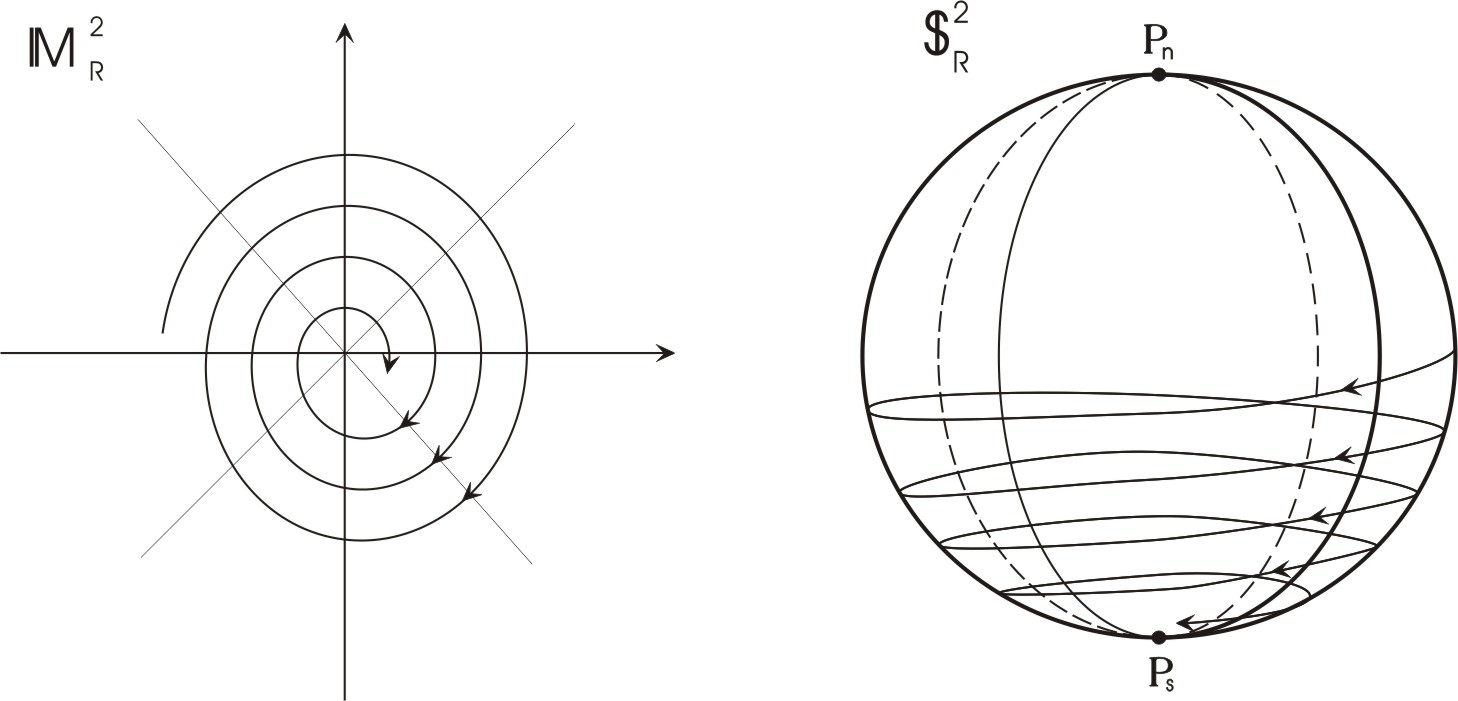}
\caption{ The M\"obius-loxodromic orbits on $\mathbb{M}_R^2$ and on the sphere.}
\label{loxodromic-action}
\end{figure}

We consider the one dimensional parameterized set of loxodromic M\"obius transformations acting in
$\mathbb{M}_R^2$ given by,
\begin{equation}\label{eq:Mobius-loxodromic}
f_{G_l(t)}(z)= \lambda^2 (t) \, z= \phi(t) \, e^{it}  \, z(t),
\end{equation}
where $\phi(t)$ is a real nonnegative function defined in a suitable interval. Figure
\ref{loxodromic-action} shows the M\"obius-loxodromic (spiral) orbits of the action of the
one-parametric set $\{G_l(t)\}$ in  $\mathbb{M}_R^2$ and the corresponding
helicoidal (loxodromic) orbits in the two dimensional sphere when $\phi(t) <1$.

\begin{Definition} A solution $\mathbf{z}(t)=(z_1(t), z_2(t), \cdots, z_n(t))$ of equation (\ref{eq:motiongral})
is called M\"obius-loxodromic solution if it is invariant under the one-parametric set of M\"obius
transformations given in (\ref{eq:Mobius-loxodromic}).
\end{Definition}

Among the whole set of M\"obius-loxodromic solutions we study two particular families
 obtained by combining two different kinds  of homothetic solutions with
those of elliptic type analyzed in section \ref{sec:hyperbolic-solutions}.
The first combined type of homothetic solutions are the hyperbolic motions obtained in section
\ref{sec:hyperbolic-solutions}, which when combine with the M\"obius-elliptic ones
arise the family of asymptotic M\"obius loxodromic motions as it can be seen in subsection
\ref{subsec:Asymptotic-loxodromic} below.
The second combined type of homothetic solutions are the so called totally geodesic, which are
defined  in subsection \ref{subsec:totally-geodesic}, and when combine with the M\"obius-elliptic
solutions arise the class of homographic M\"obius-loxodromic motions as it can be seen in subsection
\ref{subsec:Homographic-loxodromic}.

\subsection{Asymptotic M\"obius-loxodromic solutions}\label{subsec:Asymptotic-loxodromic}

We consider the one dimensional parametric subgroup of M\"obius
-loxodromic transformations
\begin{equation}\label{eq:moebius-asymptotic-initial}
f_{A(t)} (z) =e^{t}\, e^{it} \, z= e^{t(1+i)} \, z.
\end{equation}

\begin{Definition} A solution $\mathbf{z}(t)=(z_1(t), z_2(t), \cdots, z_n(t))$ of equation (\ref{eq:motiongral})
is called Asymptotic M\"obius-loxodromic if it is invariant under
the one-dimensional parametric subgroup
(\ref{eq:moebius-asymptotic-initial}).
\end{Definition}

The following result gives necessary and sufficient conditions in order to have the above kind of motions.

\begin{Theorem}\label{Theo:asympt-loxodromic} In the positively curved $n$--body problem, a necessary
and sufficient condition for having Asymptotic M\"obius-loxodromic solutions of
the equation (\ref{eq:motiongral}) is that the positions $z_k(t)$ of all the particles satisfy the
system of algebraic equations
\begin{equation} \label{eq:rationalsystem-asympt-loxodromic}
\frac{4 \,R^6 i \, (R^2-|z_k|^2) z_k }{(R^2+ |z_k|^2)^4} =
 \sum_{j=1, j \neq k}^n \frac{m_j \, (|z_j|^2+R^2)^2(R^2+\bar{z}_jz_k)(z_j-z_k)}{|z_j-z_k|^3 \, |R^2+ \bar{z}_j z_k|^3}
\end{equation}
and the velocity in each particle is given by the relation
\begin{equation}\label{eq:velocities-asympt-loxodromic}
2 \dot{z}_k =-(1+i) z_k,
\end{equation}
for $k=1,2,\cdots, n$. Moreover, all  solutions
$z_k=z_k(t)$ are forward asymptotic in an helicoidal way to the
origin of coordinates, which implies that there are not collisions
among the particles. \label{thm:existence-asympt-loxodromic}
\end{Theorem}

{\bf Proof.} If we suppose that $w_k  = e^{t(1+i)} \, z_k$ is a
solution of  equation (\ref{eq:motiongral}), then by differentiating
we obtain
\begin{eqnarray}
\dot{w}_k  &=& (\dot{z}_k + (1+i) \, z_k ) \, e^{t(1+i)} \nonumber \\
\ddot{w}_k  &=& (\ddot{z}_k + 2 (1+i) \dot{z}_k+ (1+i)^2 z_k ) \, e^{t(1+i)}, \nonumber \\
\end{eqnarray}
which when we substitute in the equation (\ref{eq:motiongral}) gives
us the relation
\begin{eqnarray}
&& m_k (\ddot{z}_k + 2 (1+i) \dot{z}_k+ (1+i)^2 z_k ) \, e^{t(1+i)}
- \frac{2 m_k e^{t(1-i)}\bar{z}_k e^{2t(1+i)} (\dot{z}_k + (1+i) \,
z_k )^2}{R^2+ |z_k|^2\, e^{t(1+i)}} \nonumber \\
&=& \frac{2}{\lambda (z_k \, e^{t(1+i)}, \bar{z}_k\, e^{t(1+i)})} \,
\frac{\partial U_R }{\partial \bar{z}_k}(z, \bar{z}) \, e^{-t(1-i)}. \nonumber \\
\end{eqnarray}

Again, by evaluating at $t=0$ we obtain the condition for the
infinitesimal generator of the vector field associated to the
solution $z_k=z_k(t)$ in one arbitrary point. This gives the
relation
\begin{equation}
 m_k (\ddot{z}_k + 2 (1+i) \dot{z}_k+ (1+i)^2 z_k )  -\frac{2 m_k \bar{z}_k (\dot{z}_k + (1+i) \,
z_k )^2}{R^2+ |z_k|^2} = \frac{2}{\lambda (z_k, \bar{z}_k )} \,
\frac{\partial U_R }{\partial \bar{z}_k}(z, \bar{z}),
\end{equation}
and using that $z_k=z_k(t)$ holds equation (\ref{eq:motiongral}),
then this last equation becomes into the one
\begin{equation}\label{eq:generator:asymtotic-loxodromic}
2 (1+i) \dot{z}_k+ (1+i)^2 z_k  =\frac{2 \bar{z}_k [ (1+i)^2 z_k^2 +
2 (1+i) \dot{z}_k  z_k ]}{R^2+ |z_k|^2},
\end{equation}
or equivalently, to the equation,
\begin{equation}\label{eq:generator:asymtotic-loxodromic-2}
(1+i)\left[1-\frac{2 |z_k|^2}{R^2+ |z_k|^2} \right] \, \left[ (1+i)
z_k + 2  \dot{z}_k \right] =0.
\end{equation}

The equation (\ref{eq:generator:asymtotic-loxodromic-2}) holds for
the infinitesimal conditions $ |z_k|=R$, which corresponds to the
geodesic equator, or for
\begin{equation}\label{eq:generator:asymtotic-loxodromic-4}
2 \,  \dot{z}_k = - (1+i)\, z_k
\end{equation}
which shows  that the solution through the point $z(0)$ is an
helicoidal curve converging asymptotically to the origin of
coordinates. Since the solutions move along different helicoidal
curves, then there are not collisions among the particles along these kind of solutions.

If we derive the infinitesimal condition given in equation
(\ref{eq:generator:asymtotic-loxodromic-4}), we obtain
\begin{equation}\label{eq:generator:asymtotic-loxodromic-5}
2 \,  \ddot{z}_k =i \, z_k
\end{equation}
which when is substituted into the equations of motion
(\ref{eq:motiongral}) gives  the system
(\ref{eq:rationalsystem-asympt-loxodromic}).
This ends the proof of the Theorem. \qed

By straightforward computations, very similar to the ones  we did in
 the above sections,  we can prove  the following result for the
asymptotic M\"obius-loxodromic  motion of two particles in $\mathbb{M}_R^2$.

\begin{Proposition} \label{theo:two-bodies-asymptotic}
Suppose that two particles with masses $m_1$ and $m_2$ in $\mathbb{M}^2_{R}$ and
positions $z_1(t)$ and $z_2(t)$
move as an asymptotic M\"obius-loxodromic solution. Then, the masses are equal iff the particles are located on opposite sides of the same circle, that is, $z_1(t)=-z_2(t)$.
\end{Proposition}

 For the curved $3$--body problem on $\mathbb{M}_R^2$, we get the following results whose proofs follow
 by straightforward computations, we omit them the proof of
both results here.

\begin{Proposition} \label{theo:lagrangian-asymptotic}
Three equal masses moving on $\mathbb{M}^2_R$ form an asymptotic M\"obius-loxodromic solution iff
the particles form an equilateral triangle for all time.
\end{Proposition}

\begin{Proposition}\label{theo:eulerian-asypmtotic} A necessary and sufficient condition in order to have a
collinear asymptotic M\"obius-loxodromic solution in $\mathbb{M}^2_{R}$, with one particle at the origin $z_3(t)=0$,
is that the masses of the other two particles with positions $z_1(t)$ and $z_2(t)$  be equal,
and that $z_1(t)=-z_2(t)$.
\end{Proposition}


\subsection{Totally geodesic solutions}\label{subsec:totally-geodesic}

In order to find another type of M\"obius-loxodromic solutions for equation (\ref{eq:motiongral}), we begin by searching
homothetic solutions of  the form $z_k(t)= \phi (t) \, z_{k,0}$ lying on the meridians
through the point  $z_{k,0}$  and parameterized as a geodesic curve. For this,
we impose that the homothetic function $\phi (t)$ in  (\ref{eq:Mobius-loxodromic}) must satisfy the relations
\begin{eqnarray}\label{eq:cond-homothetic-function}
\ddot{\phi} z_{k,0} -\frac{2  \bar{z}_{k,0} z_{k,0}^2 \, \phi \, \dot{\phi}^2}{R^2+ \phi^2 \, |z_{k,0}|^2} &=& 0, \\
\frac{\partial U_R }{\partial \bar{z}_k}\left(\phi (t) z_{k,0}, \phi (t)  \bar{z}_{k,0}\right) &=& 0, \nonumber
\end{eqnarray}
for $k=1,2, \cdots, n$.

\smallskip

In this way we obtain the function which state the geodesic parametrization of such solutions.

\begin{Lemma}\label{lemma:homothetic-function} The first second order differential equation
of system (\ref{eq:cond-homothetic-function}) can be integrated by
quadratures by a smooth real function $\phi=\phi (t)$ with initial conditions $\phi(0)=1$ and $\dot{\phi}(0)=-1$,
which when is substituted into the second equation of (\ref{eq:cond-homothetic-function}) holds the equality.
\end{Lemma}

{\bf Proof.} For $z_{k,0} \neq 0$, the first differential equation is equivalent to,
\[ \frac{\ddot{\phi}}{\dot{\phi}} = \frac{2  |z_{k,0}|^2 \, \phi \, \dot{\phi}}{R^2+ \phi^2 \, |z_{k,0}|^2},\]
which when is integrated becomes into
\begin{equation}\label{eq:first-homotetic-function}
\frac{\dot{\phi}}{R^2+  \phi^2 \, |z_{k,0}|^2}= C_2,
\end{equation}
for a suitable real constant $C_2$.
If we integrate directly the equation (\ref{eq:first-homotetic-function})
we obtain the relation
\begin{equation}\label{eq:second-homotetic-function}
\phi(t)= \frac{R^2}{|z_{k,0}|^2} \tan \left( C_2|z_{k,0}| \, t + C_1\right),
\end{equation}
for other  suitable constant of integration $C_1$. Now, using the initial conditions,
relation (\ref{eq:second-homotetic-function}) we obtain  the particular solution
\begin{equation}\label{eq:particular-homotetic-function}
\phi(t)= \frac{R^2}{|z_{k,0}|^2} \tan \left( \arctan \left(\frac{|z_{k,0}|^2}{R^2}\right)-
 \left(\frac{R^2 |z_{k,0}|}{|z_{k,0}|^4+R^4}\right) \, t \right).
\end{equation}

Since function (\ref{eq:particular-homotetic-function}) makes that all particles move locally along
meridians in a geodesic way, then along all such solutions the second equation in
(\ref{eq:cond-homothetic-function}) vanishes. This remark ends the proof of the Lemma. \qed

\begin{Definition} A homothetic solution $\mathbf{z}(t)=(z_1(t), z_2(t), \cdots, z_n(t))$ of equation (\ref{eq:motiongral})
is called  totally geodesic if  all the particles move along geodesic curves which hold equations (\ref{eq:cond-homothetic-function}).
\end{Definition}

We obtain the following result for these kind of motions.

\begin{Theorem}\label{Theo:homothetic-solutions} In the positively curved $n$--body problem, a necessary
and sufficient condition for having totally geodesic solutions of the equation (\ref{eq:motiongral})
is that all particles move along the straight lines through the origin of coordinates (meridians) and
the whole set of solutions satisfy the system of equations
\begin{eqnarray}\label{eq:cond-hyperbolic-function}
0 &=& \frac{\partial U_R }{\partial \bar{z}_k}\left(z_{k}, \bar{z}_{k}\right)  \nonumber \\
&=&
 \sum_{j=1, j \neq k}^n \frac{m_k m_j \, ( R^2+
 |z_k|^2)(|z_j|^2+R^2)^2(R^2+\bar{z}_jz_k)(z_j-z_k)}{4 R^2 \, |z_j-z_k|^3 \, |R^2+ \bar{z}_j z_k|^3} \nonumber \\
&=&  \sum_{j=1, j \neq k}^n \frac{m_j \, (|z_j|^2+R^2)^2}{4 R^2 \, |z_j-z_k|
\, |R^2+ \bar{z}_j z_k| \, (\bar{z}_j-\bar{z}_k) \, (R^2+\bar{z}_k z_j)},
\end{eqnarray}
for $k=1,2, \cdots, n$.
\end{Theorem}

{\bf Proof.} By hypothesis,  any particle with mass $m_k$ must holds the relations
(\ref{eq:cond-homothetic-function}),
which implies that the corresponding solution moves with the geodesic parametrization given in  Lemma
\ref{lemma:homothetic-function} in the direction of the
meridian thorough such point and the origin of coordinates. Since one solution of the equation (\ref{eq:motiongral})
moves along one  geodesic curve if and only if the right hand side (equation (\ref{eq:cond-hyperbolic-function})) vanishes,
the claim follows. \qed

When $n=2$, the following result shows the no-existence of totally geodesic orbits,
an unintuitive and surprising result.

\begin{Theorem}\label{Theo:homothetic-2-bodies} For the positively curved two-body problem there are not
totally geodesic orbits.
\end{Theorem}

{\bf Proof.} Let $m_1$ and $m_2$ be two arbitrary masses in the space  $\mathbb{M}^2_{R}$ moving under
equations (\ref{eq:motiongral}) and conforming an homothetic solution, that is, the particles are moving
along two (possibly different)  geodesics.
The system (\ref{eq:cond-hyperbolic-function}) becomes into,
\begin{eqnarray}\label{eq:homothetic-2-bodies}
0 &=&  \frac{m_2 (R^2+ z_1 \bar{z}_2)(z_2-z_1)}{ |z_2-z_1|^3 \, |R^2+ \bar{z}_2 z_1|^3}, \nonumber \\
0 &=&  \frac{m_1  (R^2+ z_2 \bar{z}_1)(z_1-z_2)}{|z_1-z_2|^3 \, |R^2+ \bar{z}_1 z_2|^3},
\end{eqnarray}
with the infinitesimal condition (\ref{eq:generator:hyperbolic3}).

\smallskip

We remark that avoiding collisions $z_i \neq z_j$ and antipodal points $R^2+ \bar{z}_i z_j \neq 0$,
a condition necessary and sufficient  or the existence of nontrivial solutions of the linear system
(\ref{eq:homothetic-2-bodies}) for the masses is that the determinant
of such system vanishes, that is,
\begin{eqnarray}
0 &=& \frac{(R^2+ z_1 \bar{z}_2)(z_2-z_1)}{ |z_2-z_1|^3 \, |R^2+
\bar{z}_2 z_1|^3}\frac{(R^2+ z_2 \bar{z}_1)(z_2-z_1)}{|z_1-z_2|^3 \, |R^2+ \bar{z}_1 z_2|^3}, \nonumber \\
 &=& \frac{1}{|z_2-z_1|^2 \, (\bar{z}_2-\bar{z}_1)^2 \,|R^2+ \bar{z}_2 z_1|^4} \nonumber \\
\end{eqnarray}
which never holds for $z_1, z_2 \in \mathbb{M}^2_R$. This proves the claim. \qed

In the spirit of \cite{Diac-Perez} we obtain conditions for having
totally geodesic solutions of the system (\ref{eq:motiongral}) which
have the same shape configuration for all time $t$ and all the particles move along geodesic curves.
That is, all the particles are located at any time on the same euclidian circle and they do not rotate.
Moreover, the property to be totally geodesic for a solution implies that the curved gradient vanishes along all
the curves conforming a solution.

\smallskip

Let $m_1$, $m_2$ and $m_3$ be three masses in the space  $\mathbb{M}^2_{R}$
moving along a particular totally geodesic solution, such that the  corresponding configuration
satisfy $r(t)= |z_1(t)|=|z_2(t)|$, with  $z_3(t)$ fixed at the
origin. Along the solution the following system holds (these kind of orbits are called Eulerian
solutions).
\begin{eqnarray}\label{eq:homothetic-3-eulerian-bodies}
0 &=&  \frac{m_2  (R^2+ z_1 \bar{z}_2)(z_2-z_1)}{|z_2-z_1|^3 \, |R^2+ \bar{z}_2 z_1|^3}
- \frac{ m_3 \, z_1}{ R^4 \, |z_1|^3}, \nonumber \\
0 &=&  \frac{m_1  (R^2+ z_2 \bar{z}_1)(z_1-z_2)}{ |z_2-z_1|^3 \, |R^2+ \bar{z}_2 z_1|^3}
 -  \frac{m_3  \, z_2}{R^4 \,  |z_2|^3}, \nonumber \\
0 &=&  \frac{m_1  \, z_1}{|z_1|^3 }
 +  \frac{m_2 \, z_2}{|z_2|^3},
\end{eqnarray}

\begin{Theorem}\label{theo:homothetic-3-eulerian-bodies} A necessary and sufficient condition for having an
 Eulerian totally geodesic solution of system (\ref{eq:homothetic-3-eulerian-bodies})
 is that the masses $m_1$ and $m_2$ be equal and
the corresponding particles are always located on  the same circle but in opposite sides.
\end{Theorem}

{\bf Proof.} A necessary and sufficient condition in the linear system (\ref{eq:homothetic-3-eulerian-bodies})
for having nontrivial solutions for the masses is that the principal determinant vanishes,
but since we avoid collisions ($z_1 \neq z_2$) and antipodal points ($R^2+ \bar{z}_j z_k \neq 0$),
it is equivalent to the equality
\begin{equation}\label{eq:homothetic-3-bodies-lagrangian-1}
 z_2 \bar{z}_1- z_1 \bar{z}_2=0.
\end{equation}

Without loss of generality we can put
$z_1(t)=r(t) \, e^{i\, \theta_1}$, $z_2(t)=r(t) \, e^{i\, \theta_2}$ in equation
(\ref{eq:homothetic-3-bodies-lagrangian-1}), for suitable constant
angles $\theta_1$ y $\theta_2$. Such relation is equivalent to the trigonometric equation
\begin{equation}\label{eq:homothetic-3-bodies-lagrangian-2}
\sin (\theta_2-\theta_1) = 0
\end{equation}
which, avoiding collisions holds iff $\displaystyle \theta_2 =\ \theta_1 + \pi$.

Therefore we can do $z_1(t)=r(t)$ and  $z_2(t)=-r(t)$. Substituting these values in the third equation
 of system (\ref{eq:homothetic-3-eulerian-bodies})  we obtain that
the masses $m_1$ and $m_2$ must be equal. This ends the proof of the Theorem. \qed

\begin{Remark}
Assuming that the masses $m_1$ and $m_2$ are equal, and that they are symmetrically located
on the same circle,  in \cite{Diac-Perez}, the authors obtain a homothetic Eulerian solution.
In this sense Theorem \ref{theo:homothetic-3-eulerian-bodies} is stronger, since in principle the masses
could be situate on different geodesic meridians and not necessarily symmetrically located.
\end{Remark}

{\begin{Corollary}  There is not totally geodesic solution for the restricted
eulerian positively curved problem.
\end{Corollary}

{\bf Proof.} The restricted eulerian positively curved problem is obtained when
the third mass tends to zero, $ m_3 \to 0$ in the system (\ref{eq:homothetic-3-eulerian-bodies}).
Such system becomes in the more simple two-body system
(\ref{eq:homothetic-2-bodies}) plus the equation
\[ \frac{m_1  \, z_1}{|z_1|^3 }  +  \frac{m_2 \, z_2}{|z_2|^3}=0,\]
which together with Theorem  \ref{theo:homothetic-3-eulerian-bodies} implies that $m_1=m_2$
and $z_1=-z_2$. A direct substitution of these pair of relations in any  equation
of (\ref{eq:homothetic-3-eulerian-bodies}) allows us to a contradiction. This proves
the claim and ends the proof. \qed

\smallskip

Let $m_1$, $m_2$ and $m_3$ be three masses in the space  $\mathbb{M}^2_{R}$ moving along
any  totally geodesic solution. Then the following system holds,
\begin{eqnarray}\label{eq:homothetic-3-lagrangian-bodies}
0 &=&  \frac{m_2  (R^2+ z_1 \bar{z}_2)(z_2-z_1)}{ |z_2-z_1|^3 \, |R^2+ \bar{z}_2 z_1|^3}
+  \frac{m_3 (R^2+ z_1 \bar{z}_3)(z_3-z_1)}{ |z_3-z_1|^3 \, |R^2+ \bar{z}_3 z_1|^3}, \nonumber \\
0 &=&  \frac{m_1  (R^2+ z_2 \bar{z}_1)(z_1-z_2)}{|z_1-z_2|^3 \, |R^2+ \bar{z}_1 z_2|^3}
 +  \frac{m_3  (R^2+ z_2 \bar{z}_3)(z_3-z_2)}{|z_3-z_2|^3 \, |R^2+ \bar{z}_3 z_2|^3}, \nonumber \\
0 &=&  \frac{m_1  (R^2+ z_3 \bar{z}_1)(z_1-z_3)}{|z_1-z_3|^3 \, |R^2+ \bar{z}_1 z_3|^3}
 +  \frac{m_2  (R^2+ z_3 \bar{z}_2)(z_2-z_3)}{ |z_2-z_3|^3 \, |R^2+ \bar{z}_2 z_3|^3}.
\end{eqnarray}

Suppose that a particular configuration satisfy $r(t)= |z_1(t)|=|z_2(t)|=|z_3(t)|$.
We will show that if, in general, three arbitrary masses are located on the same circle, they generate
a totally geodesic solution if the masses form an equilateral triangle for all time.

\begin{Theorem}\label{theo:homothetic-3-lagrangian-bodies} The configuration of a totally geodesic
solution of the positively curved $3$--body problem for arbitrary masses located on the same circle for all time is
always an equilateral triangle.
\end{Theorem}

{\bf Proof.} A necessary and sufficient condition in the linear system (\ref{eq:homothetic-3-lagrangian-bodies})
for having nontrivial solutions for the masses is that the principal determinant vanishes,
but since we avoid collisions and antipodal points, it is equivalent to the equality
\begin{equation}\label{eq:homothetic-3-bodies-1}
 (R^2+ z_1 \bar{z}_3)(R^2+ z_2 \bar{z}_1)(R^2+ z_3 \bar{z}_2)- (R^2+ z_1 \bar{z}_2)(R^2+ z_1 \bar{z}_3)(R^2+ z_2 \bar{z}_3)=0,
\end{equation}

Without loss of generality we can do
$z_1(t)=r(t)$ in equation (\ref{eq:homothetic-3-bodies-1}), $z_2(t)=r(t) \, e^{i\, \theta_2}$ and
$z_3(t)=r(t) \, e^{i\, \theta_3}$, for suitable constant
angles $\theta_2$ and $\theta_3$, such relation is equivalent to
the system of trigonometric equations
\begin{eqnarray}\label{eq:homothetic-3-bodies-2}
\sin \theta_2 - \sin \theta_3 + \sin (\theta_3-\theta_2) =0, \nonumber \\
\cos \theta_2 + \cos \theta_3 + \cos (\theta_3-\theta_2) =0,
\end{eqnarray}
which avoiding collisions hold iff  $\displaystyle \theta_2 =\frac{2 \pi}{3}$
and $\displaystyle \theta_3  =\frac{4 \pi}{3}$.

In another words, equation (\ref{eq:homothetic-3-bodies-1}) holds iff the configuration is
an equilateral triangle. \qed

\subsection{Homographic M\"obius-loxodromic
solutions}\label{subsec:Homographic-loxodromic}

We consider the one dimensional parametric set of M\"obius--loxodromic transformations
\begin{equation}\label{eq:moebius-homographic-initial}
f_{H(t)} (z) =\phi (t) \, e^{it} \, z,
\end{equation}
where $\phi (t)$ is the geodesic homothetic function obtained in
Lemma \ref{lemma:homothetic-function}.

\begin{Definition} A solution $\mathbf{z}(t)=(z_1(t), z_2(t), \cdots, z_n(t))$ of equation (\ref{eq:motiongral})
is called homographic M\"obius-loxodromic if it is invariant under
the one-dimensional parametric set
(\ref{eq:moebius-homographic-initial}).
\end{Definition}

The next step is to obtain algebraic conditions on the position of
the particles for having Homographic M\"obius-loxodromic solutions
in the positively curved problem, these are given in the following
result.

\begin{Theorem}\label{Theo:helicoidal-motions} Consider $n$ point particles
with masses $m_1,m_2, \cdots, m_n>0$ moving in $\mathbb{M}^2_{R}$. A
necessary and sufficient condition for  the solution
$\mathbf{z}(t)=(z_1(t), z_2(t), \cdots, z_n(t))$ of
(\ref{eq:motiongral})  to be a homographic M\"obius-loxodromic
solution is that the coordinates satisfy the following system of
rational equations.
\begin{equation} \label{eq:rational-helicoidal-system}
\frac{R^6 \, (3i-1)^2(R^2-|z_k|^2) z_k }{2 (R^2+ |z_k|^2)^4  }
= - \sum_{j=1, j \neq k}^n \frac{m_j (R^2+ |z_j|^2)^2 (R^2+z_k\bar{z}_j)(z_j-z_k)}{|z_j-z_k|^3
\, |R^2+ \bar{z}_j z_k|^3},
\end{equation}
for $k=1,2, \cdots, n$, and for the condition in the velocities $\displaystyle \dot{z}_k = \frac{3i-1}{4} \, z_k$.
In particular, if all the particles are located on the same euclidian circle, $|z_k|=|z_j|=r(t)=r$, then the system of
algebraic equations
\begin{equation} \label{eq:rational-helicoidal-system-2}
\frac{R^6 \, (3i-1)^2(R^2-r^2) z_k }{2(R^2+ r^2)^6  }
= - \sum_{j=1, j \neq k}^n \frac{m_j
(R^2+z_k\bar{z}_j)(z_j-z_k)}{|z_j-z_k|^3 \, |R^2+ \bar{z}_j z_k|^3}
\end{equation}
must hold.
\end{Theorem}

{\bf Proof.} Let us suppose that $z_k=z_k(t)$ is a solution  of (\ref{eq:motiongral}) and let $w_k$ be
the conformal  function given by
\[w_k(t)=\phi (t) \, e^{it} \, z_k(t), \]
then, by differentiating we obtain
\begin{eqnarray}
\dot{w}_k &=& (\dot{\phi} \,  z_k + i \,\phi  \,  z_k
+ \phi  \, \dot{z}_k)\, e^{it},  \\
\ddot{w}_k &=& (\ddot{\phi}  \, z_k + 2i \, \dot{\phi}  \, z_k
+ 2 \dot{\phi}\,  \dot{z}_k - \phi \, z_k + 2 i \, \phi  \dot{z}_k
+ \phi   \, \ddot{z}_k )\, e^{it}. \nonumber
\end{eqnarray}

Using the fact that  $w_k$ is a solution of (\ref{eq:motiongral}) and the relation
$\displaystyle \frac{\bar{z}_k}{\bar{w}_k}= \frac{e^{it}}{\phi} $,  substituting in  (\ref{eq:motiongral})
we obtain
\begin{eqnarray}\label{eq:helicoidal-solutions-1}
&& m_k \,(\ddot{\phi}  \, z_k + 2i \, \dot{\phi}  \, z_k
+ 2 \dot{\phi}\,  \dot{z}_k - \phi \, z_k + 2 i \, \phi  \dot{z}_k
+ \phi   \, \ddot{z}_k )\, e^{it} \nonumber \\
&-& \frac{2m_k \phi e^{-it} \bar{z}_k}{R^2+\phi^2 |z_k|^2}\left(\dot{\phi} \,  z_k + i \,\phi  \,  z_k
+ \phi  \, \dot{z}_k \right)^2 \, e^{2it}  \nonumber \\
&=& \frac{(R^2+\phi^2 |z_k|^2)^2}{2R^4} \frac{\partial U_R}{\partial \bar{z}_k} \, \frac{e^{it}}{\phi}.
\end{eqnarray}

Since $e^{it} \neq 0$, by using the condition for the first equation of (\ref{eq:cond-homothetic-function}) in equation (\ref{eq:helicoidal-solutions-1}), we get that for all time $t$ the function $\phi$ and the solution $z_k$ must
satisfy the system
\begin{eqnarray}\label{eq:helicoidal-solutions-2}
&&m_k \,( 2i \, \dot{\phi}  \, z_k+ 2 \dot{\phi}\,  \dot{z}_k - \phi \, z_k + 2 i \, \phi  \dot{z}_k
+ \phi   \, \ddot{z}_k ) \nonumber \\
&-& \frac{2m_k \phi \, \bar{z}_k}{R^2+\phi^2 |z_k|^2}
\left(-\phi^2 \, z_k^2 + \phi^2 \, \dot{z}_k ^2+   2 i \, \dot{\phi} \phi z_k ^2
 +   2 \dot{\phi} \phi z_k \, \dot{z}_k  + 2i \, \phi^2 \,z_k \, \dot{z}_k\right)      \nonumber \\
&=& \frac{(R^2+\phi^2 |z_k|^2)^2}{2R^4 \, \phi} \frac{\partial U_R}{\partial \bar{z}_k}.
\end{eqnarray}

Again, we search the conditions for having the infinitesimal generators of
 the vector field for such solutions by doing $t=0$. Introducing the initial conditions
$\phi(0)=1$ and $\dot{\phi}(0)=-1$ into equation (\ref{eq:helicoidal-solutions-2}),
we obtain the algebraic rational system of equations
\begin{eqnarray}\label{eq:helicoidal-solutions-3}
&& m_k \,(- 2i  \, z_k- 2   \dot{z}_k -  z_k + 2 i \,   \dot{z}_k
+  \ddot{z}_k ) \nonumber \\
&-& \frac{2m_k  \, \bar{z}_k}{R^2+ |z_k|^2}
\left( - z_k^2 + \dot{z}_k ^2 -   2 i \,  \phi z_k ^2
 -   2  z_k \, \dot{z}_k  + 2i \, z_k \, \dot{z}_k  \right)     \nonumber \\
&=& \frac{(R^2+ |z_k|^2)^2}{2R^4} \frac{\partial U_R}{\partial \bar{z}_k}.
\end{eqnarray}
Now, since $m_k \neq 0$ and $z_k$ is a solution for equation (\ref{eq:motiongral}) for all time $t$,
the above system (\ref{eq:helicoidal-solutions-3}) give us
\begin{equation}\label{eq:helicoidal-solutions-4a}
 (- 2i  \, z_k- 2   \dot{z}_k -  z_k + 2 i \,   \dot{z}_k )
= \frac{2 \, \bar{z}_k}{R^2+ |z_k|^2}
\left( - z_k^2 -   2 i \, z_k ^2 -   2  z_k \, \dot{z}_k  + 2i \, z_k \, \dot{z}_k  \right).
\end{equation}

Last equation (\ref{eq:helicoidal-solutions-4a}) holds if it satisfy at least one of the following conditions
\begin{eqnarray}\label{eq:helicoidal-solutions-5}
|z_k|&=& R, \nonumber \\
\dot{z}_k &=& \frac{2i+1}{2i-2} \, z_k =\frac{3i-1}{4} \, z_k .
\end{eqnarray}

We observe that the first equation in (\ref{eq:helicoidal-solutions-5}) corresponds to the geodesic
circle of radius $R$. The second is a first order differential equation
(the infinitesimal conditions for having such motions), whose integrals are
loxodromic (helicoidal) curves in $\mathbb{M}^2_{R}$, parameterized by
\begin{equation}\label{eq:helicoidal-solutions-6}
z_k(t)=z_{k,0} \, e^{(\frac{3i-1}{4}) \, t},
\end{equation}
 the same second equation in (\ref{eq:helicoidal-solutions-5})  gives us the relation of velocities that
must be hold for obtaining such solutions.

\smallskip

Finally,  deriving the second equation we obtain
\begin{equation}\label{eq:helicoidal-solutions-7}
\ddot{z}_k  = \frac{3i-1}{4} \, \dot{z}_k = \frac{(3i-1)^2}{16} \, z_k.
\end{equation}

A direct substitution of the second relation of (\ref{eq:helicoidal-solutions-5}) and
the relation (\ref{eq:helicoidal-solutions-7}) into equations of motion (\ref{eq:motiongral})
gives   the system of algebraic rational equations (\ref{eq:rational-helicoidal-system}).
This ends the proof of Theorem \label{Theo:helicoidal-motions}. \qed

\smallskip

 We remark that the second equation in (\ref{eq:helicoidal-solutions-5}) implies that any
M\"obius-loxodromic solution intersects all the geodesics rays (the meridians on the sphere)
in a constant angle. This is the reason for call those solutions loxodromic.

\medskip

Now we study some examples of M\"obius-loxodromic solution, we start with the
particular case of two particles in the space $\mathbb{M}^2_{R}$
where $r= r(t)=|\phi (t)|= |z_k(t)|$  (for  $k=1,2$). Actually this example has been widely studied
in \cite{Perez}, from where
we know  that two particles of masses $m_1$ and
$m_2$ moving on the same circle of radius $r \neq R$ form a relative
equilibria iff the masses are equal and they are located at opposite
sides of the circle. For
M\"obius-loxodromic solutions we have the following result.

\begin{Theorem} For the positively curved two-body problem
with equal masses there are not homographic M\"obius-loxodromic solutions.
\end{Theorem}

{\bf Proof.} In this case, we have that  $m_1=m_2=m$, then the system of rational equations
(\ref{eq:rational-helicoidal-system}) becomes into
\begin{eqnarray} \label{eq:rational-helicoidal-two-bodies}
\frac{R^6 (3i-1)^2(R^2-r^2) z_1 }{ (R^2+ r^2)^6 } &=&  - \frac{m
(R^2+z_1\bar{z}_2)(z_2-z_1)}{|z_2-z_1|^3 \, |R^2+ \bar{z}_2 z_1|^3} \nonumber \\
\frac{R^6 (3i-1)^2(R^2-r^2) z_2 }{ (R^2+ r^2)^6 } &=&  - \frac{m
(R^2+z_2\bar{z}_1)(z_1-z_2)}{|z_1-z_2|^3 \, |R^2+ \bar{z}_1 z_2|^3} \nonumber \\
\end{eqnarray}

 If we avoid collisions and conjugated points from the system (\ref{eq:rational-helicoidal-two-bodies}) we
 obtain the condition
\begin{equation}
-(R^2+z_1\bar{z}_2)\bar{z}_1 =(R^2+z_2\bar{z}_1)\bar{z}_2,
\end{equation}
which, since $|z_1(t)|=r(t)=r$ and $|z_2(t)|=r(t)=r$, it is equivalent to
\begin{equation}
-R^2 \, (\bar{z}_1 + \bar{z}_2) = r^2 \, (\bar{z}_2 + \bar{z}_1),
\end{equation}
which does not hold for any mass. \qed

\smallskip

We remark that this last result matches with the result of Theorem \ref{Theo:homothetic-2-bodies} which shows
that there are not totally geodesic solutions in the two body problem in $\mathbb{M}^2_{R}$.

\medskip

Now we use the results obtained in sections \ref{sec:relative-equilibria} and
\ref{sec:hyperbolic-solutions}, for obtaining M\"obius-loxodromic
(homographic) solutions for the particular case of three particles in $\mathbb{M}^2_{R}$
when $r= r(t)=|\phi (t)|= |z_k(t)|$  (for  $k=1,2,3$). We start generalizing the Eulerian orbits
of celestial mechanics.

\begin{Theorem}\label{Theo:eulerian-helicoidal} A necessary and sufficient condition for a
solution of the curved $3$-body problem in
$\mathbb{M}^2_{R}$, with one particle at the origin and the other two with equal mass located on the same
  geodesic, to be a M\"obius-loxodromic
solution is that the particles are located always at opposite sides of the same
circle of radius $r=r(t) \neq R$.
\end{Theorem}

{\bf Proof.}  The proof of this Theorem follows by straightforward computations, using the
conditions of the configuration in equations  (\ref{eq:rational-helicoidal-system}). The point to remark here
is the conditions on the velocities (\ref{eq:helicoidal-solutions-5}) of
 the particles in the same circle. This is given by $\displaystyle \dot{z}_k = \frac{3i-1}{4} \, z_k$.
 \qed

\smallskip

For the generalization of the Lagrangian solutions of celestial mechanics, we have the following result,
whose proof is also by straightforward computations, we omit it here, getting exactly the same condition
(\ref{eq:helicoidal-solutions-7}) on the velocities of the particles as in the previous Theorem.

\begin{Theorem}\label{Theo:Lagrangian-helicoidal} Let us consider a
configuration of 3 equal masses on the same circle of radius
$r=r(t)$  in $\mathbb{M}^2_R$, then a necessary and  sufficient
condition to have  a M\"obius-loxodromic solution is that the
particles form always an equilateral triangle.
\end{Theorem}

\begin{Remark}
The M\"obius-loxodromic solutions are the similar to the homographic solutions studied in \cite{Diac-Perez}, many
 results for the case $n=3$ were proved first in that paper, nevertheless we think that the proofs present in
 this work are easier, showing the big advantage to work with intrinsic coordinates and
 the use of complex variable for the computations.
\end{Remark}

\subsection*{Acknowledgments} This work has been
partially supported by CONACYT, M\'exico, Grant 128790.

\end{document}